\documentclass[10pt,twosided]{article}
\usepackage[tbtags]{amsmath}
\usepackage{amssymb}
\allowdisplaybreaks[4]

\usepackage{color}

\usepackage{amssymb,color}

\definecolor{c20}{rgb}{0.,0.7,0.}
\definecolor{c30}{rgb}{0.,0.,1.}
\definecolor{c40}{rgb}{1,0.1,0.7}
\definecolor{c50}{rgb}{1,0,0}
\definecolor{c60}{rgb}{0,0.9,0.1}

\usepackage{amssymb,color}

\newcommand{\kb}[1]{\boldsymbol{#1}}
\newcommand{\vk}[1]{\kb{#1}}
\def\kal#1{{\cal{ #1}}}

\newcommand{\abs}[1]{\lvert #1 \rvert}

\newcommand{\E}[1]{\mathbb{E}\left(#1\right)}

\newcommand{\pk}[1]{\mathbb{P} \left( #1 \right) }

\newcommand{\R}{\mathbb{R}}

\newcommand{\N}{\mathbb{N}}

\newcommand{\ldot}{,\ldots,}

\newcommand{\todis}{\stackrel{d}{\to}}

\newcommand{\BQN}{\begin{eqnarray}}
\newcommand{\EQN}{\end{eqnarray}}
\newcommand{\BQNY}{\begin{eqnarray*}}
\newcommand{\EQNY}{\end{eqnarray*}}

\newcommand{\BS}{\begin{sat}}
\newcommand{\ES}{\end{sat}}
\newcommand{\BT}{\begin{theo}}
\newcommand{\ET}{\end{theo}}
\newcommand{\BK}{\begin{korr}}
\newcommand{\EK}{\end{korr}}

\newcommand{\BD}{\begin{de}}
\newcommand{\ED}{\end{de}}

\newcommand{\BIT}{\begin{itemize}}
\newcommand{\EIT}{\end{itemize}}
\newcommand{\BDI}{\begin{description}}
\newcommand{\EDI}{\end{description}}

\newcommand{\BRM}{\begin{remarks}}
\newcommand{\ERM}{\end{remarks}}

\newcommand{\BEL}{\begin{lem}}
\newcommand{\EEL}{\end{lem}}

\newtheorem{theo}{Theorem}[section]
\newtheorem{sat}[theo]{Proposition}
\newtheorem{de}[theo]{Definition}
\newtheorem{lem}[theo]{Lemma}

\newtheorem{korr}[theo]{Corollary}

\newtheorem{remarks}[theo]{Remarks}

\newcommand{\nelem}[1]{{Lemma \ref{#1}}}

\newcommand{\netheo}[1]{{Theorem \ref{#1}}}
\newcommand{\nekorr}[1]{{Corollary \ref{#1}}}

\newcommand{\prooftheo}[1]{ \textbf{Proof of Theorem} \ref{#1}:}

\newcommand{\prooflem}[1]{\textbf{Proof of Lemma} \ref{#1}:}
\newcommand{\proofkorr}[1]{\textbf{Proof of Corollary} \ref{#1}:}

\newcommand{\COM}[1]{}

\newcommand{\QED}{\hfill $\Box$}

\def\rw{\rightarrow}

\def\IF{\infty}

\date{}

\def\oo{(1+o(1))}

\def\LT{\left}
\def\RT{\right}

\def\oo{(1+o(1))}
\def\MM{\mathcal{M}}
\def\HMM{\widehat{\MM}}
\def\Piter{\mathcal{P}}
\def\PTT{ \Upsilon_n (u) }
\def\Var{\mathrm{Var}}
\def\Corr{\mathrm{Corr}}
\def\ovX{\overline{Y}}
\def\Cov{\mathrm{Cov}}
\def\Vol{\mathrm{mes}}

\def\en{\epsilon}
\def\QQ{\mathbb{Q}}

\begin{document}

\title{\bf Extremes of Chi-square Processes with Trend}

\bigskip
\author{{Peng Liu\thanks{ School of Mathematical Sciences and LPMC, Nankai University, Tianjin 300071, China} }$^{,\dagger}$\ \ \   and\ \ Lanpeng Ji\thanks{Department of Actuarial Science,
University of Lausanne, UNIL-Dorigny 1015 Lausanne, Switzerland} 
}

 \maketitle
\vskip -0.61 cm

  \centerline{\today{}}

\bigskip
{\bf Abstract:} 
This paper studies  the supremum of chi-square processes with trend over a threshold-dependent-time horizon. Under the assumptions that the chi-square process is generated from a centered self-similar Gaussian process and the trend function is modeled by a polynomial function, we obtain the exact  tail asymptotics of the supremum of the chi-square process with trend. These results are of interest in applications in engineering, insurance, queuing  and statistics, etc. Some possible extensions of our results are also discussed.




\section{Introduction}

Let $\{Y(t), t\ge0\}$ be  a centered self-similar Gaussian process with almost surely (a.s.) continuous sample paths and index $H\in(0,1)$, i.e., $\Var(Y(t))=t^{2H}$  and for any $a>0$ and any $s,t\ge0$
\BQNY
\Cov(Y(at),Y(as))=a^{2H}\Cov(Y(t),Y(s)).
\EQNY
It has been shown that self-similar Gaussian processes such as fractional
 Brownian motion (fBm), sub-fractional Brownian motion and bi-fractional
 Brownian motion are quite useful in applications in engineering, telecommunication, insurance, queueing, finance, etc., see
 \cite{debicki2002ruin, EM02, HP99, KK04, Norros94, RT06} and the references therein.

Let $\beta,c$ be two positive constants.  In this paper we are interested in the tail asymptptics of the supremum of a chi-square process with trend given by
\BQN \label{eq:psiT}
\psi_{T}(u) =  \pk{\sup_{t\in[0,T]}\left(\sum_{i=1}^nb_i^2 Y_i^2(t)-ct^\beta\right)>u},\ \ \ u\to\IF,
\EQN
where $Y_i, i=1,\cdots,n$ are independent copies of the centered self-similar Gaussian process $Y$, and  $1=b_1=\cdots=b_k>b_{k+1}\ge b_{k+2}\ge\cdots\ge b_n>0$. Here  $T>0$ can be a finite constant, infinity, and eventually we allow $T=T_u, u>0$ to be a threshold-dependent
positive deterministic function.

One motivation for considering \eqref{eq:psiT} stems from its applications in engineering sciences, see \cite{Lindgren1980} and the references therein. More precisely, let $\vk{X}(t)=(X_1(t),\cdots,X_n(t)), t\ge0$ be a vector Gaussian load
process. Of interest is the probability of exit
\BQNY
\pk{\vk{X}(t)\not\in\vk{ S}_u(t),\ \ \text{for}\ \text{some}\ t\in[0,T]},
\EQNY
where the {\it time-dependent safety region} $\vk{S}_u(t), t\ge0$ is defined by
$$
\vk{S}_u(t)=\Biggl\{(x_1,\cdots,x_n)\in\R^n: {\sum_{i=1}^n x_i^2}\le h(t,u)\Biggr\}
$$
with $h(t,u), t, u\ge0$ some positive function.
Various models for $\vk{X}$ and $h(t,u)$ (especially, $h(t,u)\equiv u$) have been discussed in the literature (e.g., \cite{ Berman82, AL1990, PI1994b, Pit96, HAJI2014}) for the case that $T\in(0,\IF)$. In this framework, $\psi_T(u)$ corresponds to the model with $\vk{X}=(b_1Y_1,\cdots,b_nY_n)$ and  $h(t,u)=u+ct^\beta$. As one of the new features of this contribution, we shall deal with different types of $T=T_u,u\ge0$; see Section \ref{Sec_Tu}.

Another
motivation  stems from its applications in insurance. Specifically,  the surplus process of an insurance company can be modeled by
\BQN\label{eq:R}
R_u(t)=u+ct^\beta- \sum_{i=1}^nb_i^2 Y_i^2(t), \ \   t\ge0,
\EQN
where $u$ is the  initial reserve, $ct^\beta$ models the total premium received up to
time $t$, and $\sum_{i=1}^nb_i^2 Y_i^2(t)$ represents the total amount of  aggregated
claims up to time $t$ from $n$ different types of risks. In this framework, $\psi_T(u)$ is called {\it ruin probability} which is the most important measure of risk of the insurance company; see, e.g., \cite{AsmAlb10, Roletal1999}.
\COM{It is worth mentioning that in \cite{HUSC2006} and
 \cite{DESI2011} the following risk process
\BQNY
\widetilde{R}_u(t)=u+ct^\beta-\sum_{i=1}^n b_iY_{i}(t),\ \  t\ge0,
\EQNY
is chosen to model the surplus process. While from practical point of view, when modeling the aggregated claims, it is more reasonable to consider $\sum_{i=1}^nb_i^2 Y_i^2(t)$ (instead of $\sum_{i=1}^nb_i Y_i(t)$) as in \eqref{eq:R}. Moreover, allowing the time horizon to be adjusted by the initial reserve level $u$ of the portfolio (i.e., $T=T_u$) is one of the new features of this contribution. The motivation for doing so is the insurance rational that  if the company allocates a high initial capital $u$ to a specific insurance portfolio, then the time horizon that this portfolio is not ruined, say with at least $99\%$ probability, should be closely related to the level $u$.}
 Note that the model in  \eqref{eq:R} is also related with   the framework of fluid queue; see, e.g., \cite{DESI2011}.

Finally, we remark that the study of $\psi_T(u)$  also gives some insight into the study of some limiting test statistics. In \cite{DUWE2014}, it is shown that a test statistic converges weakly to
\BQN\label{M2}
\sup_{t\in(0,1)}\left(\frac{U(t)^2}{2t(1-t)}-C(t)-\upsilon D(t)\right),
\EQN
where $\{U(t),t\in[0,1]\}$ is a standard Brownian bridge, $C(t)=\ln \left(1-\ln (1-(2t-1)^2)\right)$, $D(t)=\ln (1+C(t)^2)$ and $\upsilon >1.$
Apparently, the above process involved is a chi-square process with trend.   Asymptotical results for the tail probability of \eqref{M2} is very interesting from statistical point of view; see, e.g., \cite{JAPI2011}. See also \cite{KHFK14} and the references therein for more applications of chi-type processes in statistics.


Outline of the rest of the paper:   Section \ref{Sec_Pre} is concerned about some preliminary results. In \netheo{ThmCS} we show the tail asymptotics of the supremum of a chi-square process generated from a non-stationary Gaussian process which extends some results in  \cite{PI1994b, HAJI2014}; \nelem{GF} derives a  Fernique-type inequality for certain Gaussian random fields.
In Section \ref{SecIF} we concentrate on the asymptotics of \eqref{eq:psiT} over an infinite-time horizon (i.e., $T=\IF$). Under a local stationary condition on the correlation of
the self-similar process $Y$ (see  \eqref{eq:locstaXH}), in \netheo{ThmInf} we derive the
asymptotics of $\psi_\IF(u)$.
 Section \ref{Sec_Tu} is devoted to the symptotics of \eqref{eq:psiT} over a threshold-dependent-time horizon (i.e., $T=T_u$ a positive deterministic function). 
 As a corollary, we also obtain  approximations of the conditional first passage time of the process defined in \eqref{eq:R}.
 Finally, in Section \ref{Sec_ExDis} possible extensions of our results are discussed. We show that general results can also be obtained for the model  where $Y_i$'s are  independent but not necessarily identical and for the model with a more general correlation structure (for $Y$) than that in  \eqref{eq:locstaXH}.

\section{ Preliminaries} \label{Sec_Pre}

Let $\{X(t),t\ge 0\}$ be a centered non-stationary Gaussian process with a.s. continuous sample paths. In the following, unless otherwise stated, $T$ is considered to be a positive finite constant. We impose the following typical  assumptions on the Gaussian process $X$ (see \cite{Pit96}):

{\bf Assumption I}:  The standard deviation function $\sigma_X(\cdot):=\sqrt{\Var(X(\cdot))}$ of $X$ attains
its maximum (assumed to be 1)
over $[0,T]$ at the unique point $t=t_0\in[0,T]$.  Further, there exist some positive constants
$\mu,a$ such that
\BQNY\label{eq:VarX}
\sigma_X(t)= 1 -a\abs{t-t_0}^{\mu}\oo, \quad t\to t_0.
\EQNY
{\bf Assumption II}: There exist some $\nu\in (0,2], d>0$ such that
\BQNY\label{eq:rX}
r_X(s,t)=\Corr(X(s), X(t))=1- d|t-s|^{\nu}\oo, \quad  s,t  \to t_0.
\EQNY
{\bf Assumption III}: There exist some positive constants $G$, $\gamma$ and $\rho$ such that
\BQNY\label{eq:IncX}
\E{(X(t)-X(s))^{2}}  \leq   G|t-s|^{\gamma}
\EQNY
holds for all $s,t\in[t_0-\rho,t_0+\rho]\cap [0,T]$.

For such a centered non-stationary Gaussian process $X$, it is known that  (see, e.g., \cite{PitPri78},
Theorem D.3 in \cite{Pit96}  or Theorem 2.1 in \cite{DHJ14})

\BQN\label{PitThm}
\pk{\sup_{t\in[0,T]}X(t)>u}=\MM_{\nu,\mu,d,a}  \frac{1}{\sqrt{2\pi}}u^{\LT(\frac{2}{\nu}-\frac{2}{\mu}\RT)_+ -1}\exp\LT(-\frac{u^2}{2}\RT)\oo, \  \ u \rw \IF,
\EQN
where $(x)_+=\max(0,x)$, and, with $I_{(\cdot)}$ denoting the indicator function,
\BQN\label{MM}
\MM_{\nu,\mu,d,a}= \left\{
          \begin{array}{ll}
  d^{1/\nu} a^{-1/\mu} \Gamma(1/\mu+1)(1+I_{(t_0\not\in \{0,T\})}) \mathcal{H}_{\nu}, & \hbox{if } \nu<\mu,\\
 \Piter_{\nu}^{\frac{a}{d}}, &  \hbox{if } \nu=\mu,\\
1, &  \hbox{if } \nu>\mu.
              \end{array}
            \right.
\EQN
Here $ \mathcal{H}_{\nu} \in (0,\IF)$ is the {\it Pickands constant} defined by
\BQNY\label{pick}
\mathcal{H}_{\nu}=\lim_{S\rightarrow\infty}\frac{1}{S}\E{ \exp\biggl(\sup_{t\in[0,S]}\Bigl(\sqrt{2}B_\nu(t)-t^{\nu}\Bigr)\biggr)}
\EQNY
with $\{B_\nu(t),t\in\R\}$ a standard  fBm  defined on $\R$ with  Hurst index  $\nu/2\in (0,1]$; and $\Piter_{\nu}^{\frac{a}{d}}\in(0,\IF)$ is the
  {\it Piterbarg  constant} defined by
  \BQN\label{PitC}
  \Piter_{\nu}^{\frac{a}{d}}=\widehat{\Piter}_{\nu}^{\frac{a}{d}}I_{(t_0\in(0,T))}+\widetilde{\Piter}_{\nu}^{\frac{a}{d}}I_{(t_0\in\{0,T\})}\in(0,\IF)
\EQN
  with
\BQNY
&&\widehat{\Piter}_{\nu}^{\lambda}=\lim_{S_1, S_2\rw\IF} \Piter_{\nu}^{\lambda}[-S_1,S_2],\ \   
 \widetilde{\Piter}_{\nu}^{\lambda}= \lim_{S\rw\IF} \Piter_{\nu}^{\lambda}[0,S] =\lim_{S\rw\IF} \Piter_{\nu}^{\lambda}[-S,0], \\
&&  \Piter_{\nu}^{\lambda}[-S_1,S_2] =\E{\exp\left(\sup_{t\in[-S_1,S_2]}\left(\sqrt{2}B_{\nu}(t)-(1+\lambda)\abs{t}^{\nu}\right)\right)}, \ \ \ \lambda>0, \max(S_1,S_2)>0.
\EQNY

 We refer to \cite{Pit96, debicki2002ruin, DeMan03, DikerY} for the properties and generalizations of the Pickands-Piterbarg and related constants.

Let  $\{\chi_{n,\vk{b}}^2(t),t\ge0\}$ be a
chi-square process with $n$ degrees of freedom defined by
\BQN \label{eq:Gchi}
\chi_{n,\vk{b}}^2(t)\ =\  \sum_{i=1}^n b_i^2 X_i^2(t),\quad  t\ge0,
\EQN
where $b_i>0$,  $1\le i\le n$ and $\{X_i(t), t\ge0\}$, $1\le i\le n$, are independent copies of the centered Gaussian process $ X$ satisfying assumptions {\bf I--III}.
As an analogue of \eqref{PitThm}, \cite{HAJI2014} derived the following tail asymptotics for $\chi_{n,\vk{1}}^2$:
\BQN\label{AlPi}
\pk{\sup_{t\in[0,T]} \chi_{n,\vk{1}}^2(t) >u} =\MM_{\nu,\mu,d,a} u^{\LT(\frac{1}{\nu}-\frac{1}{\mu}\RT)_+ } \PTT\oo,\ \ \ u\to\IF,
\EQN
where
$$
 \PTT:=\pk{ \chi_{n,\vk{1}}^2(0) >u}= \frac{2^{(2-n)/2}}{\Gamma(n/2)}u^{n/2-1}\exp\left(-\frac{u}{2}\right),\ \ u\ge0.
$$
The result in \eqref{AlPi} was derived by using a similar double-sum method as in \cite{PI1994b}. As shown in \cite{PI1994b, HAJI2014} the usage of the double-sum method for the chi-square process is usually technical, since we have to deal with the supremum of a Gaussian random field with variance function attaining its maximum on an infinite set; see also \cite{DHJ14c} for a recent result in this direction.
Below, we present a general result on the tail asymptotics of $\chi_{n,\vk{b}}^2$ allowing for different $b_i$'s. 
The next result  may not be surprising (see \cite{PI1994b, HAJI2014}), but it turns out that the proof is far from trivial. 
As we will see  the following result is crucial when dealing with the tail asymptotics of the supremum of the chi-square process with trend; two other extensions of \netheo{ThmCS} will be discussed in Section \ref{Sec_ExDis}.
\BT\label{ThmCS}
Let  $\{\chi_{n,\vk{b}}^2(t),t\ge0\}$ be a
chi-square process defined as above with generic $X$ satisfying assumptions {\bf I--III}. If $1=b_1=\cdots=b_k>b_{k+1}\ge b_{k+2}\ge\cdots\ge b_n>0$, then, as $u\to\IF$,
\BQN
\pk{\sup_{t\in[0,T]} \chi_{n,\vk{b}}^2(t) >u}
 = \prod_{i=k+1}^n (1-b_i^2)^{-1/2}\MM_{\nu,\mu,d,a} u^{\LT(\frac{1}{\nu}-\frac{1}{\mu}\RT)_+ } \Upsilon_k(u)\oo.\label{CSt02}
\EQN
\ET

We conclude this section with a  Fernique-type inequality, which will be used in the proof of our main result. 
The proof of it is quite similar to the classical Fernique's inequality (see, e.g., \cite{leadbetter1983extremes}). We refer to  \cite{MWX} for new developments on the  Fernique-type inequality.

\BEL\label{GF}
Let $\{\xi(\vk{t}), \vk{t} \in[0,1]^n\}$ be a centered Gaussian process with a.s. continuous sample paths and $\Var(\xi(0))=\sigma^2\geq 0$. Suppose that
\BQN\label{APP1}
\E{\left(\xi(\vk{t})-\xi(\vk{s})\right)^2}\leq \mathbb{Q}\sum_{i=1}^n|t_i-s_i|^{\alpha_i}
\EQN
holds for all $\vk{t},\vk{s}\in[0,1]^n$, with some  constants $\QQ>0, \alpha_i>0, 1\leq i\leq n$. Then, for all $x>0$ 
\BQNY
\mathbb{P}\left(\sup_{\vk{t}\in[0,1]^n}\xi(\vk{t})>x\right)\leq 2^{n+1}\exp\LT( -\frac{c^* x^2}{\QQ}\RT)+2^{-1}\exp\LT(-\frac{x^2}{8{\sigma}^2}\RT),
\EQNY
where $c^* = \Bigl( 2 n \sum_{p=0}^\IF\left((p+1)2^{-(p+1)\min_{1\leq i\leq n}\alpha_i +1}\right)^{1/2}\Bigr)^{-2}$,
and if $\sigma^2=0$ then the second term on the right-hand side disappears.
\EEL

\section{Infinite-time Horizon}\label{SecIF}
In this section we shall focus on the asymptotics of
\BQN \label{eq:psiIF}
\psi_{\IF}(u) = \pk{\sup_{t\in[0,\IF)}\sum_{i=1}^nb_i^2 Y_i^2(t)-ct^\beta>u},\ \ \ u\to\IF,
\EQN
with $Y_i$'s are the centered self-similar Gaussian processes as discussed in Section 1.
Throughout the paper, for technical reasons we assume that $\beta>2H$.
As demonstrated  in \cite{HP99, HP08} it is useful to define, for $\beta>2H$ and  $c>0$
\BQN\label{eq:X}
Z_i(t)=\frac{Y_i(t)}{\sqrt{1+ct^\beta}},\ \  t\ge0,\ \ \ 1\le i\le n.
\EQN
Indeed, by self-similarity of $Y_i$'s, for any $u>0$
\BQN\label{eq:ss}
\psi_{\IF}(u)=
\pk{\sup_{t\ge0} \sum_{i=1}^nb_i^2 Z_i^2(t) >u^{1-\frac{2H}{\beta}}}.
\EQN
Let
$\sigma_Z(t)= \sqrt{\Var(Z_1(t))}$ . It is noted   that
$\sigma_Z(t)$
attains its maximum on $[0,\IF)$ at the unique point
\BQNY
t_0=\LT(\frac{2H}{c(\beta-2H)}\RT)^{\frac{1}{\beta}}
\EQNY
and
\BQN\label{VAR}
\sigma_Z(t)=A^{-1/2}\LT(1- \frac{B }{4A}(t-t_0)^2 \oo\RT),\ \  t\to t_0
\EQN
with
\BQN\label{eq:AB}
  A=\left(\frac{2H}{c(\beta-2H)}\right)^{-2H/\beta}\frac{\beta}{\beta-2H}, \ \ \ B=2\left(\frac{2H}{c(\beta-2H)}\right)^{-2(H+1)/\beta}H\beta.
\EQN

 In the rest of the paper we assume  {\it  local stationarity}
for the standardized Gaussian process $\ovX(t):=Y(t)/t^H, t>0$ in a
neighborhood of the point $t_0$, i.e.,
\BQN\label{eq:locstaXH}
\lim_{s\to t_0,t\to t_0}\frac{ \E{(\ovX(s)-\ovX(t))^2}}{ \abs{s-t}^\alpha}=Q>0
\EQN
holds for some  $\alpha \in(0,2)$.
Condition \eqref{eq:locstaXH} is common in the literature; most of the known self-similar Gaussian processes (such as fBm, sub-fBm, and bi-fBm) satisfy \eqref{eq:locstaXH}, see e.g., \cite{HAJI2014}. Note that the local stationarity at $t_0$ and  the  self-similarity of the process $Y$ imply the local stationarity
at any point $r\in(0,\IF)$. 
\COM{$$
\lim_{s\to r, t\to r}\frac{ \E{(\ovX(s)-\ovX(t))^2}}{ \abs{s-t}^\alpha}=\LT(\frac{t_0}{r}\RT)^\alpha Q.
$$}

Next we present our main result concerning the tail asymptotics of the supremum of the self-similar chi-square process with trend over an infinite-time horizon.

\BT\label{ThmInf} Suppose that the generic process $\{Y(t), t\geq0\}$ is a centered  self-similar Gaussian process with index $H\in(0,1)$ and correlation function satisfying (\ref{eq:locstaXH}). If $\beta>2H$, then
\BQNY
\psi_{\IF}(u)&=&2^{1-1/\alpha} Q ^{1/\alpha}A^{1/\alpha}B^{-1/2}\pi^{1/2}\mathcal{H}_\alpha\prod_{i=k+1}^n (1-b_i^2)^{-1/2}\\
&&\times u^{(1-2H/\beta)\LT( {1}/{\alpha}- {1}/{2}\RT)} \Upsilon_k(Au^{1-2H/\beta})(1+o(1)),\ \ \ u\rw\IF.
\EQNY
\ET



\section{Threshold-dependent-time Horizon}\label{Sec_Tu}

In this section we are concerned about the asymptotics of
\BQNY
\psi_{T_u}(u) = \pk{\sup_{t\in[0,T_u]}\sum_{i=1}^nb_i^2 Y_i^2(t)-ct^\beta>u},\ \ \ u\to\IF.
\EQNY
Throughout this section we shall adopt the same notation as in Section \ref{SecIF}. In addition, define 
 $$
 B(u)=2^{1/2}B^{-1/2}u^{\frac{H+1}{\beta}-\frac{1}{2}},\ \ \ u>0.
 $$
  In what follows, the following two scenarios of $T_u>0$ will be discussed:
\begin{itemize}
\item[i)] The short time horizon: $\lim_{u\rw\IF}\frac{T_u}{u^{1/\beta}}=s_0\in[0,t_0)$;
\item[ii)]  The long time horizon: $\lim_{u\rw\IF}\frac{T_u-t_0u^{1/\beta}}{B(u)}=x\in(-\IF,\IF$].
\end{itemize}
Clearly, $T=\IF$ is included in scenario ii) and $T\in(0,\IF)$ is covered by scenario i). 

We present below our main result of this section.

\BT\label{THRP}
  Suppose that the generic   process $\{Y(t),t\geq 0\}$ is a centered self-similar Gaussian process with index $H\in(0,1)$ and correlation function satisfying (\ref{eq:locstaXH}). Assume further that $\beta>2H$. We have, as $u\to\IF$,

i) If $\lim_{u\rw\IF}\frac{T_u}{u^{1/\beta}}=s_0\in[0,t_0)$, then
\BQNY
\psi_{T_u}(u) = \prod_{i=k+1}^n (1-b_i^2)^{-1/2}\MM_{\alpha,1,\frac{Q}{2}t_0^\alpha,D} \left(\frac{u+cT_u^\beta}{T_u^{2H}}\right)^{\LT(\frac{1}{\alpha}-1\RT)_+ } \Upsilon_k\LT(\frac{u+cT_u^\beta}{T_u^{2H}}\RT)(1+o(1)), 
\EQNY
where the constant $\MM_{\alpha,1,\frac{Q}{2}t_0^\alpha,D}$ is given as in \eqref{MM} with $D=\frac{2H-c(\beta-2H)s_0^\beta}{2(1+cs_0^\beta)}$.\\
ii) If  $\lim_{u\rw\IF}\frac{T_u-t_0u^{1/\beta}}{B(u)}=x\in(-\IF,\IF]$, then
\BQNY
\psi_{T_u}(u)=\psi_\IF(u)\Phi(x)(1+o(1)), 
\EQNY
  where the asymptotics of $\psi_\IF(u)$ is given in \netheo{ThmInf} and $\Phi(\cdot)$ denotes the standard normal distribution function.
\ET

As a corollary of \netheo{THRP} we  derive  an approximation of the first passage time of the chi-square process with trend, which goes in line with
e.g., \cite{HP08, DHJ13a, HJ14c}. Precisely,
define
$$
\tau_u=\inf\{t\geq 0: R_u(t)\leq 0\}\ \ (\text{with}  \inf\{\varnothing\}=\IF )
$$
to be the first passage time to 0 of the process $\{R_u(t), t\ge0\}$ defined in \eqref{eq:R}. Denote by $\todis$ convergence in distribution when the argument tends to infinity, and let $E$ be a unit mean exponential random variable and $\mathcal{N}$ be a standard normal random variable. We have:
\COM{
Next we give the asymptotic distribution  of the scaled  first passage time $\tau_u$ given that the process ever passages 0 before time $T_u$. This  result  goes in line with,
e.g., \cite{HP08, DHJ13a, HJ14c},
where the approximation of the conditional  ruin time is considered.}
\BK\label{THRT}
Under the conditions and notation of Theorem \ref{THRP}  

 i) If $\lim_{u\rw\IF}\frac{T_u}{u^{1/\beta}}=s_0\in[0,t_0)$, then
\BQNY
\frac{(2H-\frac{c\beta s_0^\beta}{1+cs_0^\beta})(u+cT_u^\beta)}{2T_u^{2H+1}}(T_u-\tau_u)\Bigl\lvert (\tau_u\leq T_u)&\todis& E,\ \ \ u\rw\IF.
\EQNY
ii) If  $\lim_{u\rw\IF}\frac{T_u-t_0u^{1/\beta}}{B(u)}=x\in(-\IF,\IF]$, then
\BQNY
\frac{\tau_u-t_0u^{1/\beta}}{B(u)}\Bigl\lvert (\tau_u\leq T_u)&\todis& \mathcal{N}\Bigl\lvert (\mathcal{N}\leq x),\ \ \ u\rw\IF.
\EQNY
\EK

\section{Extensions  \& Discussions}\label{Sec_ExDis}

In Section 3 and Section 4, we have derived asymptotical results for the case where the chi-square process is generated from a self-similar Gaussian process. %
In this section, we shall discuss two possible extensions: (a) instead of independent copies of a self-similar Gaussian process we shall consider  
independent but non-identical self-similar Gaussian processes; (b)  instead of polynomial function $\abs{t-s}^\alpha$ in   \eqref{eq:locstaXH} we consider a regularly varying function $K^2(\abs{t-s})$ with index $\alpha\in(0,2]$.

As we have seen, Theorem  \ref{ThmCS} and Theorem \ref{Thmxx} are fundamental for the proofs of our results in the last two sections. Asymptotical results for the extended chi-square processes (as in the cases (a) and (b))  with trend will follow similarly if  corresponding extended  results for Theorems \ref{ThmCS} and \ref{Thmxx} are available. Therefore, it is sufficient at this point to present only an extension of Theorem \ref{ThmCS}; corresponding extension for \netheo{Thmxx} can also be obtained.

\subsection{Non-identical Gaussian processes $X_i$'s}
Let $\{X_i(t), t\geq 0\}, 1\leq i\leq k$ be independent copies of the  a.s. continuous  Gaussian process $X$ satisfying assumptions {\bf I--III} with the parameters therein, and let $\{X_i(t), t\geq 0 \},  k+1\leq i\leq n$ be independent copies of another  a.s. continuous Gaussian process $X^{(1)}$ satisfying assumption {\bf III} with  parameter $\gamma_1$ instead of $\gamma$. Moreover, we suppose that the standard deviation function $\sigma_{X^{(1)}}(\cdot)$ attains its maximum   1  over $[0,T]$ at  $t_0$ as well. Besides,  $\{X_i(t), t\geq 0\}, 1\leq i\leq k,$  and $\{X_i(t), t\geq 0 \},  k+1\leq i\leq n$ are assumed to be independent. Define also
\BQNY
\chi_{n,\vk{b}}^2(t)  =   \sum_{i=1}^n b_i^2 X_i^2(t),\quad  t\ge0,
\EQNY
with $1=b_1=\cdots=b_k\geq b_{k+1}\geq \cdots \ge b_n >0$.

\BT\label{THEX}
Let $\{\chi_{n,\vk{b}}^2(t),t\geq 0\}$ be a chi-square process defined as above. 
If $\gamma\geq \nu$ and $\gamma_1\geq \nu$, then we have, as $u\rw\IF$,
\BQN
\pk{\sup_{t\in[0,T]} \chi_{n,\vk{b}}^2(t) >u}
 = \prod_{i=k+1}^n (1-b_i^2)^{-1/2}\MM_{\nu,\mu,d,a} u^{\LT(\frac{1}{\nu}-\frac{1}{\mu}\RT)_+ } \Upsilon_k(u)\oo.
\EQN
\ET

\begin{remarks}
a) Suppose that the generic processes $X$ and $X^{(1)}$ are both fBm with indexes $H\in(0,1)$ and $H_1\in(0,1)$, respectively. If $H_1\geq H$, then the conditions of the last theorem are fulfilled. 

b) From the proof of the last theorem we can see  the assumption that $\{X_i(t), t\geq 0 \},  k+1\leq i\leq n$ are  identical (in distribution) is not really necessary; here to simplify the notation we chose to work under this assumption.
\end{remarks}

\subsection{General correlation structure}
First, we formulate the general assumption about the correlation structure of the generic  Gaussian process $X$.   \\
{\bf Assumption II'}: There exists some $K(\cdot)$,  a regularly varying function at 0 with index $\nu/2\in(0,1]$,  such that
\BQNY\label{eq:rXK}
r_X(s,t)=\Corr(X(s), X(t))=1-  K^2(|t-s|)(1+o(1)), \quad  s,t  \to t_0.
\EQNY
Next, we introduce some further notation.
Let
$q(u)=\overleftarrow{K}(u^{-1/2})$ be the inverse function of $K(\cdot)$ at point $u^{-1/2}$ (assumed to exist asymptotically). It follows that $q(u)$ is a regularly function at infinity with index $-1/\nu$ which can be further expressed as $q(u)=u^{-1/\nu}L(u^{-1/2})$, with
$L(\cdot)$  a slowly varying function at $0$. According  to the values of $L(u^{-1/2})$ as $u\to\IF$, we consider the following three scenarios:\\ 
\textbf{C1:} $\mu>\nu$, or $\mu=\nu$ and $\lim_{u\rw\IF} L(u^{-1/2})=0$;\\
\textbf{C2:} $\mu=\nu$ and $\lim_{u\rw\IF }L(u^{-1/2})=\mathcal{L}\in(0,\IF)$;\\
\textbf{C3:} $\mu<\nu$, or $\mu=\nu$ and $\lim_{u\rw\IF}L(u^{-1/2})=\IF.$

We present below our second extension of \netheo{ThmCS}.
\BT
Under the assumptions and conditions of Theorem \ref{ThmCS} with assumption {\bf II} replaced by assumption {\bf II'}, we have, as $u\rw\IF$,
\BQN
\pk{\sup_{t\in[0,T]} \chi_{n,\vk{b}}^2(t) >u}
 = \prod_{i=k+1}^n (1-b_i^2)^{-1/2}\widetilde{\MM}_{\nu,\mu,1,a}(u) u^{\LT(\frac{1}{\nu}-\frac{1}{\mu}\RT)_+ } \Upsilon_k(u)\oo,
\EQN
where
\BQNY
\widetilde{\MM}_{\nu,\mu,1,a}(u)= \left\{
          \begin{array}{ll}
   a^{-1/\mu} \Gamma(1/\mu+1)(1+I_{(t_0\not\in \{0,T\})}) \mathcal{H}_{\nu} \overleftarrow{L} (u^{-1/2}), & \hbox{for} \  \textbf{C1} ,\\
 \Piter_{\nu}^{  a\mathcal{L}^\nu }, &  \hbox{for} \  \textbf{C2},\\
1, &  \hbox{for} \  \textbf{C3}.
              \end{array}
            \right.
\EQNY
\ET
The proof of the last theorem follows by similar arguments as in the proof of Theorem \ref{ThmCS}, and thus we only give some remarks. Actually, it relies on a corresponding extension of \nelem{GP}, 
which can be done as in the proof of Theorem 8.2 in \cite{Pit96}.  Note that the difference from the classical results in \cite{Pit96} is that for the case $\mu=\nu$ three sub-cases should be considered differently (depending on the property of $L(\cdot)$). This is not observed in the study of some other Gaussian random fields, e.g., \cite{QW73} and \cite{dieker2005extremes}, where it is shown that  the substitution of a polynomial function $d\abs{t-s}^\nu$ by a regularly varying function $K^2(\abs{t-s})$ in the correlation structure of the Gaussian random fields does not influence much on the asymptotics. However, it seems not  surprising to have these sub-cases if one examines the proof of  Theorem 8.2 in \cite{Pit96}.

\section{Further Results \& Proofs}

This section is devoted to   the proofs of Theorems \ref{ThmCS}, \ref{ThmInf}, \ref{THRP} and \ref{THEX} and Corollary \ref{THRT}. 
Let in the following $\mathbb{Q}, \mathbb{Q}_i, i=1,2,...$ denote positive constants whose values may change from line to line.

First,  we present a result concerning  the tail asymptotics of the supremum of a chi-square process over a threshold-dependent time interval, which turns out to be crucial for the proofs of \netheo{ThmCS}, \netheo{THRP} and \nekorr{THRT}. The technical proof of it is deferred to Appendix.

\BT\label{Thmxx}
Let  $\{\chi_{n,\vk{b}}^2(t),t\ge0\}$ be a
chi-square process given as in  \eqref{eq:Gchi}  with generic $X$ satisfying assumptions {\bf I--II}, and $1=b_1=\cdots=b_k>b_{k+1}\ge b_{k+2}\ge\cdots\ge b_n>0$.
Let further $\Delta_x(u)=[t_0-x_1(u)u^{-2/\mu}, t_0+x_2(u)u^{-2/\mu}]$ with functions $x_i(u), i=1,2$ such that
$$\lim_{u\to\IF}x_i(u)=x_i\in[-\IF,\IF],\ \ \ \lim_{u\to\IF}x_i(u)u^{-1/\mu}=0, \ \ i=1,2.$$
If  $-x_1<x_2$, then
\BQN
\pk{\sup_{t\in \Delta_x(u)} \chi_{n,\vk{b}}^2(t) >u^2}
 = \prod_{i=k+1}^n (1-b_i^2)^{-1/2} \widehat{\MM}_{\nu,\mu,d,a}(x_1,x_2) u^{\LT(\frac{2}{\nu}-\frac{2}{\mu}\RT)_+ } \Upsilon_k(u^2)\oo\label{CSxu}
\EQN
 as $u\to\IF$, where
 \BQN\label{HMM}
\widehat{\MM}_{\nu,\mu,d,a}(x_1,x_2)= \left\{
          \begin{array}{ll}
  d^{1/\nu} a^{-1/\mu}  \mathcal{H}_{\nu} \left(G_\mu(a^{1/\mu}x_2)-G_\mu(-a^{1/\mu}x_1)\right), & \hbox{if } \nu<\mu,\\
 \Piter_{\nu}^{\frac{a}{d}}[-d^{1/\nu}x_1,d^{1/\nu}x_2], &  \hbox{if } \nu=\mu,\\
1, &  \hbox{if } \nu>\mu,
              \end{array}
            \right.
\EQN
with $G_\mu(x)=\int_{-\IF}^x e^{-|t|^\mu}dt, x>0$ for any $\mu>0$.
\ET

\prooftheo{ThmCS} Without lose of generality we shall only consider the case that $t_0\in(0,T)$. 
As in the proof of \netheo{Thmxx},
 we consider the Gaussian random field
$${Y_{\vk{b}}}(t,\vk{v})=\sum_{i=1}^nb_iX_i(t)v_i \quad$$ 
defined on $\mathcal{G}_T=[0, T]\times \mathcal{S}_{n-1}$, where $\mathcal{S}_{n-1}$ stands for the $(n-1)$-dimensional unit sphere.
We refer to (\ref{PV})--(\ref{PHC}) below for some important properties of the Gaussian random field $Y_{\vk{b}}$.
It follows that
$$\mathbb{P}\left(\sup_{t\in [0,T]}\chi^2_{n,\vk{b}}(t)>u^2\right)=\mathbb{P}\left(\sup_{(t,\vk{v})\in \mathcal{G}_T}Y_{\vk{b}}(t,\vk{v})>u\right).$$
Therefore, we shall focus on the tail asymptotics of  
$$
\mathbb{P}\left(\sup_{(t,\vk{v})\in \mathcal{G}_T}Y_{\vk{b}}(t,\vk{v})>\ u \right),\ \ \ \ u\to\IF.
$$
Next define $\Delta_u=[t_0- (\ln u/u)^{2/\mu}, t_0+ (\ln u/u)^{2/\mu}]$, $C_{u}=\{\vk{v}\in \mathcal{S}_{n-1}: v_i\in[- \ln u/u ,  \ln u/u], k+1\leq i \leq n\}$ 
and let
$$\pi_1(u):=\mathbb{P}\left(\sup_{(t,\vk{v})\in\Delta_u \times C_u}Y_{\vk{b}}(t,\vk{v})>u\right).$$
We have, for any $u>0$ and any small  $\rho>0$
\BQN \label{eq:pi1}
\pi_1(u) \leq  \mathbb{P}\left(\sup_{(t,\vk{v})\in\mathcal{G}_T}Y_{\vk{b}}(t,\vk{v})>u\right)
&\le& \pi_1(u)+\mathbb{P}\left(\sup_{(t,\vk{v})\in \mathcal{G}_T/( [t_0-\rho,t_0+\rho]  \times C_u)}Y_{\vk{b}}(t,\vk{v})>u\right)\nonumber\\
&&+\mathbb{P}\left(\sup_{(t,\vk{v})\in ([t_0-\rho,t_0+\rho]  / \Delta_u)\times C_u)}Y_{\vk{b}}(t,\vk{v})>u\right).
\EQN
Further, in view of \netheo{Thmxx}
\BQN\label{eq:chixu}
\pi_1(u)=\prod_{i=k+1}^n (1-b_i^2)^{-1/2} \widehat{\MM}_{\nu,\mu,d,a}(-\IF,\IF) u^{\LT(\frac{2}{\nu}-\frac{2}{\mu}\RT)_+ } \Upsilon_k(u^2)\oo
\EQN
as $u\to\IF.$
By \eqref{PV} and the Borell-TIS inequality (see, e.g., \cite{AdlerTaylor})
\BQN 
\mathbb{P}\left(\sup_{(t,\vk{v})\in \mathcal{G}_T/([t_0-\rho,t_0+\rho]\times C_u)}Y_{\vk{b}}(t,\vk{v})>u\right)\le \QQ\exp\LT(-\frac{(u-\mathbb{Q}_1)^2}{2(1-\delta_0)}\RT)
\EQN
holds for all $u$ large, with some constants $ \QQ>0, \mathbb{Q}_1>0$ and $\delta_0\in(0,1)$. Further, in the light of \eqref{PV}, (\ref{PHC}) and the Piterbarg inequality given in Theorem 8.1 in \cite{Pit2001} 
\BQN  \label{eq:Pit2}
\mathbb{P}\left(\sup_{(t,\vk{v})\in ([t_0-\rho,t_0+\rho] / \Delta_u)\times C_u)}Y_{\vk{b}}(t,\vk{v})>u\right)\leq \mathbb{Q}_2u^{\frac{2(n+1)}{\gamma\wedge 2}-1}\exp\left(-\frac{u^2}{2(1-(\ln u/u )^2\mathbb{Q}_3)^2}\right)
\EQN
holds for all $u$ large, with some positive constants $\mathbb{Q}_2, \mathbb{Q}_3$. Consequently, the claim for the case that $t_0\in(0,T)$ follows from \eqref{eq:pi1}--\eqref{eq:Pit2}. This completes the proof.
\QED

\prooftheo{ThmInf}
Let $T> t_0$ be some fixed large enough integer, and let
\BQNY
 \pi(u)=\mathbb{P}\left(\sup_{t\in[0,T]}\sum_{i=1}^nb_i^2Z_i^2(t)>u^{1-2H/\beta}\right),\ \
 \pi_1(u)=\mathbb{P}\left(\sup_{t\in[T,\IF)}\sum_{i=1}^nb_i^2Z_i^2(t)>u^{1-2H/\beta}\right).
\EQNY
Clearly
\BQNY
\pi(u)\leq\psi_{\IF}(u)\leq \pi(u)+\pi_1(u).
\EQNY
By the definition of $Z_i$'s we have that there exist  some constants $\QQ>0,\rho\ge0$ such that
$$
\E{(Z_1(t)-Z_1(s))^2}\le \QQ \abs{t-s}^{\alpha}
$$
holds for any $t,s\in [t_0-\rho,t_0+\rho]$. Thus, in view of
 (\ref{VAR}), (\ref{eq:locstaXH}) and Theorem \ref{ThmCS} we conclude that
\BQNY
\pi(u)=\prod_{i=k+1}^n (1-b_i^2)^{-1/2}\MM_{\alpha,2,\frac{Q}{2},\frac{B}{4A}} \left(A(u) \right)^{ \frac{2}{\alpha}-1} \Upsilon_k((A(u))^2)(1+o(1)),\ \ \ u\rw\IF,
\EQNY
where $A(u)=A^{1/2}u^{1/2-\frac{H}{\beta}}.$ Therefore, to complete the proof  it is sufficient to show that 
\BQNY
\pi_1(u)=o(\pi(u)),\ \ \ u\to\IF.
\EQNY
To this end, let $\widetilde{Y_{\vk{b}}}(t,\vk{v})=\sum_{i=1}^nb_iA^{1/2}Z_i(t)v_i, (t,\vk{v})\in[T,\IF)\times [-1,1]^n$. We have
\BQNY
\pi_1(u)=\mathbb{P}\left(\sup_{(t,\vk{v})\in[T,\IF)\times S_{n-1}}\widetilde{Y_{\vk{b}}}(t,\vk{v})>A(u)\right).
\EQNY
We split the interval $[T,\IF)$ into subintervals $[k, k+1), k\geq T.$ For every $k\ge T$, we have (set $Y_b^{*}(t,\vk{v})=\frac{\sqrt{1+ct^\beta}}{t^H}\widetilde{Y_b}(t,\vk{v})$)
\BQN\label{3EQ1}
&&\mathbb{P}\left(\sup_{(t,v\vk{v})\in[k,k+1)\times \mathcal{S}_{n-1}}\widetilde{Y_b}(t,\vk{v})>A(u)\right)\nonumber\\
&&\leq\mathbb{P}\left(\sup_{(t,\vk{v})\in[k,k+1)\times[-1,1]^n}Y_b^{*}(t,\vk{v})>\frac{\sqrt{1+c(k-1)^\beta}}{(k-1)^H}A(u)\right).
\EQN
In addition, there exists some global constant $\QQ$ such that for any $k\ge T$
\BQNY
&&\E{Y_b^{*}(t,\vk{v})-Y_b^{*}(t',\vk{v}')}^2\\
&&\leq2An\mathbb{E}\left(\overline{Y}(t)-\overline{Y}(t')\right)^2+2A
\sum_{i=1}^n(v_i-v_i^{'})^2\\
&&\leq\QQ\left(|t-t'|^\alpha+
\sum_{i=1}^n(v_i-v_i^{'})^2\right)
\EQNY
holds for all $t,t'\in[k,k+1), \vk{v},\vk{v}'\in [-1,1]^n.$
Next we split $[-1,1]^n$ into $2^n$ subsets of the form $\prod_{i=1}^n\Delta_{i}^{j_i}, j_i=1,2$, where $\Delta_{i}^{1}=[-1,0]$ and $\Delta_{i}^{2}=[0,1]$. By using \nelem{GF} we derive, for $k\ge T$ 
\BQNY
\mathbb{P}\left(\sup_{(t,\vk{v})\in[k,k+1)\times\prod_{i=1}^n\Delta_{i}^{j_i}}Y_b^{*}(t,\vk{v})>\frac{\sqrt{1+c(k-1)^\beta}}{(k-1)^H}A(u)\right)\leq 2^{n+2}e^{-\QQ_1\frac{1+c(k-1)^\beta}
{(k-1)^{2H}}(A(u))^2},
\EQNY
with $\QQ_1=\min(\frac{c^*}{\QQ},\frac{1}{8A})$. This together with (\ref{3EQ1}) yields that
\BQNY
\mathbb{P}\left(\sup_{(t,\vk{v})\in[k,k+1)\times S_{n-1}}\widetilde{Y_b}(t,\vk{v})>A(u)\right)\leq 2^{2n+2}e^{-\QQ_1\frac{1+c(k-1)^\beta}
{(k-1)^{2H}}(A(u))^2}.
\EQNY
Consequently, since $T$ was chosen large enough
\BQNY
\pi_1(u)&\leq&\sum_{k=T}^\IF 2^{2n+2}e^{-\QQ_1\frac{1+c(k-1)^\beta}{(k-1)^{2H}}(A(u))^2}\\
&\leq&2^{2n+2}\int_{T-2}^\IF e^{-\QQ_2(A(u))^2y^{\beta-2H}}dy\\
&\leq&\mathbb{Q}_3(A(u))^{-2}e^{-\QQ_2(T-2)^{\beta-2H}(A(u))^2}=o(\pi(u)) 
\EQNY
as $u\to\IF,$ where $\mathbb{Q}_3$ is a constant depending on $T$ and $\QQ_2=c\QQ_1$.
This completes the proof.
\QED

\prooftheo{THRP} Case i). We introduce a deterministic function $m(u)=\frac{u+cT_u^\beta}{T_u^{2H}}, u>0$, and  centered Gaussian processes
 $$
 W_{u,i}(t)=\frac{Y_i(t)}{\sqrt{1-\frac{c_u}{1+c_u}(1-t^\beta)}},\ \ t\ge0,\ \ \ 1\le i \le n,
 $$
with $c_u= {cT_u^\beta}/{u}, u>0$ such that $\lim_{u\to\IF}c_u=cs_0^\beta=: c_0$.
By the self-similarity of $ Y $  we have
\BQNY
\psi_{T_u}(u) = \pk{\sup_{t\in[0,1]}\sum_{i=1}^nb_i^2 W_{u,i}^2(t)>m(u)}.
\EQNY
Let further $c_0^{\pm\epsilon}=\max(c_0\pm\epsilon,0)$ and define
$$
W_i^{\pm\epsilon}(t)=\frac{Y_i(t)}{\sqrt{1-\frac{c_0^{\pm\epsilon}}{1+c_0^{\pm\epsilon}}(1-t^\beta)}},\ \ t\ge0, \ \ \ 1\le i \le n,
$$
for any sufficiently small $\epsilon>0$. Thus we have, for $u$ large enough
\BQNY
\pk{\sup_{t\in[0,1]}\sum_{i=1}^nb_i^2\left( W_{i}^{-\epsilon}(t)\right)^2>m(u)}\leq\psi_{T_u}(u)\leq\pk{\sup_{t\in[0,1]}\sum_{i=1}^nb_i^2\left( W_{i}^{+\epsilon}(t)\right)^2>m(u)}.
\EQNY
Next we consider the upper bound of $\psi_{T_u}(u)$. It follows that $\sigma_{W_1^{+\epsilon}}(t)$ attains its maximum over $[0,1]$ at the unique point $t_0=1$ and further 
 \BQNY
 \sigma_{W_1^{+\epsilon}}(t)=1-\frac{2H-(\beta-2H)c_0^{+\epsilon}}{2(1+c_0^{+\epsilon})}|t-1|(1+o(1)),\ \ \ t\rw 1,\\
 \Corr(W_1^{+\epsilon}(t),W_1^{+\epsilon}(s))=1-\frac{t_0^\alpha Q}{2}|t-s|^\alpha(1+o(1)),\ \ \ s,t\rw 1.
 \EQNY
 In addition, there exists some $\QQ>0$ such that
 \BQNY
 \mathbb{E}\left((W_1^{+\epsilon}(t)-W_1^{+\epsilon}(s))^2\right)\leq \mathbb{Q}|t-s|^{\alpha }
 \EQNY
holds for all $s,t\in[1/2,1]$. Therefore, in view of Theorem \ref{ThmCS}
 \BQNY
 \pk{\sup_{t\in[0,1]}\sum_{i=1}^nb_i^2\left( W_{i}^{+\epsilon}(t)\right)^2>m(u)}=\prod_{i=k+1}^n (1-b_i^2)^{-1/2}\MM_{\alpha,1,\frac{Q}{2}t_0^\alpha,D_{+\epsilon}} \left(m(u)\right)^{\LT(\frac{1}{\alpha}-1\RT)_+ } \Upsilon_k(m(u))(1+o(1))
 \EQNY
 as $u\to\IF$, with $D_{+\epsilon}=\frac{2H-(\beta-2H)c_0^{+\epsilon}}{2(1+c_0^{+\epsilon})}.$ Similar arguments give the same lower bound as above (with $+\en$ replaced by $-\en$) for $\psi_{T_u}(u)$, and thus
 letting $\epsilon\rw 0$ the claim in i) follows.

Case ii). Again, using the self-similarity we derive
 \BQNY
 \psi_{T_u}(u)=\pk{\sup_{t\in[0,T_uu^{-1/\beta}]} \sum_{i=1}^nb_i^2 AZ_i^2(t) >(A(u))^2},
 \EQNY
 where $A(u)=A^{1/2}u^{1/2- {H}/{\beta}}$. Let $t_u=t_0-u^{-1/2+H/\beta}\ln u$, and define
 $$\pi_{t_u}(u)=\pk{\sup_{t\in[0,t_u]} \sum_{i=1}^nb_i^2 AZ_i^2(t) >(A(u))^2}.$$
Clearly,
 \BQNY
\pk{\sup_{t\in[t_u,T_uu^{-1/\beta}]} \sum_{i=1}^nb_i^2 AZ_i^2(t) >(A(u))^2}\leq\psi_{T_u}(u)\leq\pk{\sup_{t\in[t_u,T_uu^{-1/\beta}]} \sum_{i=1}^nb_i^2 AZ_i^2(t) >(A(u))^2}+\pi_{t_u}(u).
 \EQNY
 In the following, we shall first derive the  asymptotics of the common term  on both sides of the above formula,  which will give the exact asymptotics of $\psi_{T_u}(u)$. Then we show that $\pi_{t_u}(u)$ is asymptotically negligible.
In view of (\ref{VAR}), (\ref{eq:locstaXH}) and Theorem \ref{Thmxx}, we have, for any $x\in (-\IF,\IF)$
 \BQNY
 \pk{\sup_{t\in[t_u,T_uu^{-1/\beta}]} \sum_{i=1}^nb_i^2 AZ_i^2(t) >(A(u))^2}=\psi_{\IF}(u)\Phi(x)(1+o(1)),\ \ \ u\rw\IF.
 \EQNY
Next we show that the last formula is also valid for $x=\IF$. Since, for any fixed $y\ge0$,
 \BQNY
 \pk{\sup_{t\in[t_u,t_0+yB(u)]} \sum_{i=1}^nb_i^2 AZ_i^2(t) >(A(u))^2}\leq \pk{\sup_{t\in[t_u,T_uu^{-1/\beta}]} \sum_{i=1}^nb_i^2 AZ_i^2(t) >(A(u))^2}\leq \psi_{\IF}(u),
 \EQNY
we obtain from Theorem \ref{Thmxx} that
 \BQNY
\Phi(y)&\leq&\underline{\lim}_{u\rw\IF}\frac{ \pk{\sup_{t\in[t_u,T_uu^{-1/\beta}]} \sum_{i=1}^nb_i^2 AZ_i^2(t) >(A(u))^2}}{ \psi_{\IF}(u)}\\
&\leq&\overline{\lim}_{u\rw\IF}\frac{ \pk{\sup_{t\in[t_u,T_uu^{-1/\beta}]} \sum_{i=1}^nb_i^2 AZ_i^2(t) >(A(u))^2}}{ \psi_{\IF}(u)}\leq 1.
 \EQNY
Therefore, letting $y\rw \IF$ we  conclude that
 \BQNY
 \pk{\sup_{t\in[t_u,T_uu^{-1/\beta}]} \sum_{i=1}^nb_i^2 AZ_i^2(t) >(A(u))^2}=\psi_\IF(u)(1+o(1)),\ \ \ u\rw\IF.
 \EQNY
To complete the proof  we prove that $\pi_{t_u}(u)=o(\psi_{\IF}(u))$ as $u\rw\IF.$ We have
 \BQNY
\pk{\sup_{t\in[0,t_u]} \sum_{i=1}^nb_i^2 AZ_i^2(t) >(A(u))^2}=\pk{\sup_{(t,\vk{v})\in[0,t_u]\times \mathcal{S}_{n-1}} \widetilde{Y_{\vk{b}}}(t,\vk{v}) >A(u)},
 \EQNY
 where $\widetilde{Y_{\vk{b}}}(t,\vk{v})=\sum_{i=1}^nb_iA^{1/2}Z_i(t)v_i, (t,\vk{v})\in [0, t_0+1]\times \mathcal{S}_{n-1}$. Further, there exist  some constants $\delta\in(0,1), \QQ>0$ such that
 \BQNY
 \mathbb{E}\left((\widetilde{Y_{\vk{b}}}(t,\vk{v}))^2\right)&\leq& 1-\delta<1,\ \ \ t\in[0,t_0-\rho], \vk{v}\in \mathcal{S}_{n-1},\\
 \mathbb{E}\left((\widetilde{Y_{\vk{b}}}(t,\vk{v})-\widetilde{Y_{\vk{b}}}(t,\vk{v}))^2\right)&\leq& \mathbb{Q}(|t-s|^{\alpha}+\sum_{i=1}^n(v_i-v_i^{'})^2),\ \ \ t\in[t_0-\rho, t_0+\rho], \vk{v}\in \mathcal{S}_{n-1}
 \EQNY
 hold. Therefore, as in the proof of \netheo{ThmCS}, by the Borell-TIS inequality  we have
\BQN\label{BOR}
\pk{\sup_{t\in[0,t_0-\rho]} \sum_{i=1}^nb_i^2 AZ_i^2(t) >(A(u))^2}\leq e^{-\frac{(A(u)-\mathbb{Q}_0)^2}{2(1-\delta)^2}}=o(\psi_{\IF}(u)),\ \ \ u\rw\IF,
\EQN
with $\mathbb{Q}_0=\mathbb{E}\left(\sup_{(t,\vk{v})\in[0,t_0-\rho]\times\mathcal{S}_{n-1} } \widetilde{Y_{\vk{b}}}(t,\vk{v})\right)<\IF$,
and by the Piterbarg inequality and (\ref{VAR}) (or by a direct application of \cite{HT13}, Proposition 3.2)  we have
\BQN\label{PIT}
\pk{\sup_{t\in[t_0-\rho, t_u]} \sum_{i=1}^nb_i^2 AZ_i^2(t) >(A(u))^2}\leq \mathbb{Q}_1(A(u))^{2(n+1)/\alpha}\Psi\left(\frac{A(u)}{1-\mathbb{Q}_2(A(u)^{-1}\ln A(u))^2}\right)=o(\psi_{\IF}(u))
\EQN
as $u\rw\IF$, where $\mathbb{Q}_1$ and $\mathbb{Q}_2$ are two positive constants. Consequently, we conclude from
 (\ref{BOR}) and (\ref{PIT}) that
 $$\pi_{t_u}(u)=o(\psi_{\IF}(u)),\ \ \ \ u\rw\IF.$$
  This completes the proof. \QED

 \proofkorr{THRT} Case i). For notational simplicity, we let
 $$ f(u) = \frac{2T_u^{2H+1}}{(2H-\frac{c\beta s_0^\beta}{1+cs_0^\beta})(u+cT_u^\beta)},\ \ \ u>0.
 $$
 By definition, for any $x>0$
 \BQNY
\mathbb{P}\left(\frac{T_u-\tau_u}{f(u)}>x\Bigl\lvert \tau_u\leq T_u\right)=\frac{\psi_{T_u-xf(u)}(u)}{\psi_{T_u}(u)}.
 \EQNY
Further, it follows from Theorem \ref{THRP} that
 \BQNY
 \lim_{u\rw\IF}\frac{\psi_{T_u-xf(u)}(u)}{\psi_{T_u}(u)}=\lim_{u\rw\IF}e^{\frac{u+cT_u^\beta}{2T_u^{2H}}-\frac{u+c(T_u-xf(u))^\beta}{2(T_u-xf(u))^{2H}}}
 =e^{-x},
 \EQNY
 establishing the claim in i).

Case ii).  Similarly as above, in the light of Theorem \ref{THRP} we have, for any $y\le x$
 \BQNY
 \lim_{u\rw\IF} \mathbb{P}\left(\frac{\tau_u-t_0u^{1/\beta}}{B(u)}<y\Bigl\lvert \tau_u\leq T_u\right)=\lim_{u\rw\IF}\frac{\psi_{t_0u^{1/\beta}+yB(u)}(u)}{\psi_{T_u}(u)}=\frac{\Phi(y)}{\Phi(x)}.
 \EQNY
 Thus, the proof is complete.\QED

\prooftheo{THEX} One approach is to follow a similar proof as Theorem \ref{ThmCS} by using the double-sum method. 
Here, we give another proof based on the ideas and results in \cite{HAJI2014}, \cite{PI1994b} and \cite{Pakes04}. 
We first show that
\BQN
\pk{\sup_{t\in[0,T]} \chi_{n,\vk{b}}^2(t) >u}   =  \pk{\sup_{t\in[0,T]} \chi_{k,\vk{1}}^2(t)+\sum_{i=k+1}^n b_i^2 X_i^2(t_0) >u}\oo\label{CSt0}
\EQN
holds as $u\to\IF$, which in view of Lemma 2.1 in \cite{Pakes04} is sufficient. Indeed, letting $G(u)=\pk{\sup_{t\in[0,T]} \chi_{k,\vk{1}}^2(t) \le u}$ we have from \eqref{AlPi} that
\BQNY
\lim_{u\to\IF}\frac{1-G(u+y)}{1-G(u)}=\exp\LT(-\frac{1}{2} y \RT), \ \ \forall y\in\R.
\EQNY
Further let $H(u)=\pk{\sum_{i=k+1}^n b_i^2 X_i^2(t_0) \le u}$. It is known (cf. Example 2 in \cite{KPH13}) that
$$
1-H(u)=O\LT(u^r \exp\LT(-\frac{u}{2b_{k+1}} \RT)\RT)  =o(1-G(u))
$$
for some $r\in\N$. Moreover, by choosing some $\theta\in(1/2, 1/(2b_{k+1}))$ we have that
\BQNY
\int_{0}^\IF e^{\theta x}dH(x)<\IF.
\EQNY
Therefore, by  Lemma 2.1 in \cite{Pakes04} the claim in \eqref{CSt02} follows from \eqref{CSt0}. \\
It remains to show \eqref{CSt0}. To this end, we introduce the following two Gaussian random fields:
\BQNY
{Y_{\vk{b}}}(t,\vk{v})=\sum_{i=1}^n b_i v_i X_i(t), \ \ {Z_{\vk{b}}}(t,\vk{v})=\sum_{i=1}^k   v_i X_i(t) +\sum_{i=k+1}^n b_iv_iX_i(t_0),\ \ \ t\ge0,\ \vk{v}\in\R^n.
\EQNY
\COM{It is known that (see \cite{PI1994b} or  \cite{Pit96})
$$
\underset{t\in[0,T]}\sup \chi_{n,\vk{b}}(t)=\underset{(t,\vk{v})\in\mathcal{G}_{T}}\sup {Y_{\vk{b}}}(t,\vk{v}),\ \ \underset{t\in[0,T]}\sup \sqrt{\chi_{k,\vk{1}}^2(t)+\sum_{i=k+1}^n b_i^2 X_i^2(t_0)}
=\underset{(t,\vk{v})\in\mathcal{G}_{T}}\sup {Z_{\vk{b}}}(t,\vk{v})
$$
where $\mathcal{G}_{T}=[0,T]\times\mathcal{S}_{n-1}$, with $\mathcal{S}_{n-1}$ being the unit sphere (with respect to $L_2$-norm) in $\R^n$.
}
As in the proof of \netheo{ThmCS} it is sufficient that
\BQN\label{Eq:YZ}
\pk{ \underset{(t,\vk{v})\in\mathcal{G}_{T}}\sup {Y_{\vk{b}}}(t,\vk{v}) >u}
=\pk{\underset{(t,\vk{v})\in\mathcal{G}_{T}}\sup {Z_{\vk{b}}}(t,\vk{v})>u}\oo
\EQN
holds as $u\to \IF$. Next we have that the standard deviations  $\sigma_{Y_{\vk{b}}} (t,\vk{v})$ and $\sigma_{Z_{\vk{b}}}(t,\vk{v})$ attain their absolute maximum
 (equal to 1) over $\mathcal{G}_{T}$ at all points of $C_0$ given as 
$$
C_0=\{t_0\}\times\{\vk{v}\in\mathcal{S}_{n-1}:  v_1^2+\cdots+v_k^2=1\}\subset \mathcal{G}_{T}.
$$
Further we consider the expansions of the standard deviations and the correlations of the Gaussian random fields ${Y_{\vk{b}}}$ and ${Z_{\vk{b}}}$ around the sphere $C_0$. By direct calculations we have
\BQN \label{eq:Var}
&&\sigma_{Y_{\vk{b}}} (t,\vk{v})=1-a\abs{t-t_0}^\mu\oo-\frac{1}{2}\sum_{i=k+1}^n(1-b_i^2)v_i^2\oo,\nonumber\\
&&\sigma_{Z_{\vk{b}}} (t,\vk{v})=1-a\abs{t-t_0}^\mu\oo-\frac{1}{2}\sum_{i=k+1}^n(1-b_i^2)v_i^2\oo
\EQN
hold as $t\to t_0$ and $\sum_{i=k+1}^n v_i^2\to 0$. Further, since $\gamma>\nu, \gamma_1\ge \nu$,
\BQN \label{eq:Corr}
&&r_{Y_{\vk{b}}}(t,\vk{v},s, \vk{u})=\Corr({Y_{\vk{b}}} (t,\vk{v}),{Y_{\vk{b}}} (s,\vk{u}))=1-d\abs{t-s}^\nu\oo-\frac{1}{2}\sum_{i=1}^nb_i^2(v_i-u_i)^2\oo,\nonumber\\
&&r_{Z_{\vk{b}}}(t,\vk{v},s, \vk{u})=\Corr({Z_{\vk{b}}} (t,\vk{v}),{Y_{\vk{b}}} (s,\vk{u}))=1-d\abs{t-s}^\nu\oo-\frac{1}{2}\sum_{i=1}^nb_i^2(v_i-u_i)^2\oo
\EQN
hold as $s,t\to t_0$, $\sum_{i=k+1}^n v_i^2\to 0$ and $\sum_{i=k+1}^n u_i^2\to 0$. The technical proof of (\ref{eq:Corr}) is relegated to the Appendix.

 Define a neighborhood  $C_u$ of $C_0$ as
\BQNY
C_u=\{(t,\vk{v}): a\abs{t-t_0}^\mu+\frac{1}{2}\sum_{i=k+1}^n(1-b_i^2)v_i^2<\ln u/u\}\cap \mathcal{G}_{T}.
\EQNY
By an application of Borell inequality and Piterbarg inequality as in the proof of Lemma 8.1 in \cite{Pit96} we can show that
\BQN
&&\pk{ \underset{(t,\vk{v})\in\mathcal{G}_{T}}\sup {Y_{\vk{b}}}(t,\vk{v}) >u}=\pk{ \underset{(t,\vk{v})\in C_u }\sup {Y_{\vk{b}}}(t,\vk{v}) >u}\oo,\label{eq:CuY}\\
&&\pk{ \underset{(t,\vk{v})\in\mathcal{G}_{T}}\sup {Z_{\vk{b}}}(t,\vk{v}) >u}=\pk{ \underset{(t,\vk{v})\in C_u }\sup {Z_{\vk{b}}}(t,\vk{v}) >u}\oo\label{eq:CuZ}
\EQN
hold as $u\to\IF$.
Moreover, since we are concerned about the asymptotic results it follows that 
the expansions of the standard deviations and the correlations of the Gaussian random field ${Y_{\vk{b}}}$ (or ${Z_{\vk{b}}}$) around the sphere $C_0$ are the only necessary properties influencing the asymptotics of \eqref{eq:CuY} (or \eqref{eq:CuZ}); this is due to the fact that ${Y_{\vk{b}}}$ and ${Z_{\vk{b}}}$ are Gaussian and $C_u\to C_0$ as $u\to\IF$.
Therefore, it follows from \eqref{eq:Var} and \eqref{eq:Corr} that \eqref{Eq:YZ} is established. This completes the proof. \QED



\section{Appendix}

This section is devoted to  the proofs of \netheo{Thmxx} 
and Eq. \eqref{eq:Corr}. 


We first present  a lemma  concerning the tail asymptotics of the supremum of a Gaussian random field over a threshold-dependent-time interval, which is crucial for the proof of \netheo{Thmxx}.
\BEL\label{GP}
Let $\{X(\vk{v}), \vk{v}\in   \R^n\}$ be a centered stationary Gaussian random field  with a.s. continuous sample paths and covariance function $r(\cdot)$ satisfying  
\BQNY
r(\vk{v})=1-d_1 \abs{v_1}^\alpha (1+o(1))-\sum_{i=2}^nd_iv_i^2(1+o(1)), \quad v_1^2+v_{2}^2+\cdots+v_n^2\rw 0,
\EQNY
with  $\alpha \in(0,2]$ and $ d_i>0,i\le n$.
Define
 $$X^*(\vk{v})= \frac{X(\vk{v})}{(1+c_1|v_1|^\beta)\LT(1+\sum_{i=k+1}^nc_iv_i^2\RT)},\ \ \vk{v}\in \R^n$$
  with some  $1 < k< n$ and $c_1>0, c_i>0, k+1\le i\le n,  \beta>0$.
Let $A\subset \R^{k-1}$ be a Jordan measurable set  with positive Lebesgue measure
$\Vol(A)$.
Let further $\widetilde{\Delta}_x(u)=[-x_1(u)u^{-2/\beta}, x_2(u)u^{-2/\beta}]$ with  functions $x_i(u), i=1,2$ such that
$$\lim_{u\to\IF}x_i(u)=x_i\in[-\IF,\IF],\ \ \ \lim_{u\to\IF}x_i(u)u^{-1/\beta}=0, \ \ i=1,2.$$
Denote $D(u)= \widetilde{\Delta}_x(u)\times A\times[-  \ln u/u), \ln u/u ]^{n-k}$.
If $-x_1<x_2$, then
\BQNY
\pk{\sup_{\vk{v}\in D(u)} X^*(\vk{v}) >u}
&= & 2^{-1/2}\pi^{-k/2}\prod_{i=2}^kd_i^{1/2}
\prod_{i=k+1}^n\sqrt{1+d_i/c_i }  \Vol (A)\\
&&\times \widehat{\MM}_{\alpha,\beta,d_1,c_1}(x_1,x_2)\ u^{k-2+(2/\alpha-2/\beta)_+}e^{-\frac{u^2}{2}}\oo
\EQNY
as $u\to \IF$, with $\widehat{\MM}_{\alpha,\beta,d_1,c_1}(x_1,x_2)$ given as in \eqref{HMM}.

\EEL
{\bf Proof:} Note that
\BQNY
  \mathcal{H}_{2}=\frac{1}{\sqrt{\pi}}, \ \ \widehat{\Piter}_2^{c/d}=\sqrt{1+d/c}.
\EQNY
The proof follows by  a little modification of the proof of Theorem 8.2 in \cite{Pit96}; see also Lemma 6 in \cite{PI1994b} or Theorem 3.2 in \cite{HAJI2014}.\QED

\def\eE#1{#1}
\def\mmu{\ln u/u}
\def\tU{\Theta(u)}
\def\tUU{\eE{\widetilde{\Theta}}(u)}
\def\QQ{\mathbb{Q}}
\def\NN{\mathcal{N}}

\prooftheo{Thmxx} 
The key idea is to work with
Gaussian random fields instead of analyzing chi-square processes. The first step is standard (see, e.g., \cite{Pit96}). We consider the Gaussian random field
$${Y_{\vk{b}}}(t,\vk{v})=\sum_{i=1}^nb_iX_i(t)v_i, \quad$$ 
defined on $\mathcal{G}_x(u)=\Delta_x(u)\times \mathcal{S}_{n-1}$, where $\mathcal{S}_{n-1}$ stands for the $(n-1)$-dimensional unit sphere.
Since from \cite{Pit96}
$$\sup_{t\in \Delta_x(u)} \chi_{n,\vk{b}}(t)=\sup_{(t,\vk{v})\in \mathcal{G}_x(u)}Y_{\vk{b}}(t,\vk{v})$$
we have that
$$\mathbb{P}\left(\sup_{t\in \Delta_x(u)}\chi^2_{n,\vk{b}}(t)>u^2\right)=\mathbb{P}\left(\sup_{(t,\vk{v})\in \mathcal{G}_x(u)}Y_{\vk{b}}(t,\vk{v})> u\right).$$
It follows that the standard deviation $\sigma_{Y_{\vk{b}}}$ of $Y_{\vk{b}}$ 
 attains its maximum  (equal to $1$) over $\mathcal{G}_x(u)$ only at points on $\{(t_0,\vk{v}), \vk{v}\in \mathcal{S}_{n-1}, v_i=0, k+1\leq i \leq n\}$. Furthermore, following the arguments as in \cite{PI1994b} we conclude that
 $\sigma_{Y_{\vk{b}}}$ and the correlation function $r_{Y_{\vk{b}}}$ of $Y_{\vk{b}}$  have the following asymptotic expansions:
\BQN\label{PV}
\sigma_{Y_{\vk{b}}}(t,\vk{v})=1-a|t-t_0|^\mu\oo-\sum_{i=k+1}^n\frac{1-b_i^2}{2 }v_i^2\oo 
\EQN
as $t\rw t_0$ and $v_{k+1}^2+\cdots+v_n^2\to 0$, and
\BQN\label{PCorr}
r_{Y_{\vk{b}}}(t,\vk{v},t',\vk{v}')=1-d\abs{t-t'}^\nu\oo-\sum_{i=1}^n\frac{b_i ^2}{2}(v_i-v'_i)^2\oo  
\EQN
as $t, t'\rw t_0,$ $v_{k+1}^2+\cdots+v_n^2\to 0$, and ${v'}_{k+1}^2+\cdots+{v'}_n^2\to 0$.

In addition, there exist $\delta>0, \mathbb{Q}>0$ such that
\BQN\label{PHC}
\E{\left(Y_{\vk{b}}(t,\vk{v})-Y_{\vk{b}}(t',\vk{v}')\right)^2}\leq \mathbb{Q}(\abs{t-t'}^\gamma+\sum_{i=1}^n(v_i-v_i')^2) 
\EQN
holds for all $(t,\vk{v})\in([t_0-\rho,t_0+\rho]\cap[0,T])\times \mathcal{S}_{n-1}.$
Next 
define $C_{u}:=\{\vk{v}\in \mathcal{S}_{n-1}: v_i\in[- \ln u/u ,  \ln u/u], k+1\leq i \leq n\}$.
Let
$$\pi(u):=\mathbb{P}\left(\sup_{(t,\vk{v})\in\Delta_x(u) \times C_u}Y_{\vk{b}}(t,\vk{v})>u\right).$$
We have, for any $u>0$
\BQNY
\pi(u) \leq  \mathbb{P}\left(\sup_{(t,\vk{v})\in\mathcal{G}_x(u)}Y_{\vk{b}}(t,\vk{v})>u\right)
 \le  \pi(u)+\mathbb{P}\left(\sup_{(t,\vk{v})\in \Delta_x(u)\times (\mathcal{S}_{n-1}/ C_u)}Y_{\vk{b}}(t,\vk{v})>u\right).
\EQNY
\COM{By \eqref{PV} and the Borell-TIS inequality
\BQN \label{eq:Bor}
\mathbb{P}\left(\sup_{(t,\vk{v})\in \mathcal{G}_T/(([t_0-\rho,t_0+\rho]\cap[0,T])\times C_u)}Y_{\vk{b}}(t,\vk{v})>u\right)\le \QQ\exp\LT(-\frac{(u-a)^2}{2(1-\delta_0)}\RT)
\EQN
holds for all $u$ large, with some constants $a>0, \QQ>0$ and $\delta_0\in(0,1)$.
 }
 Further, in the light of \eqref{PV}, (\ref{PHC}) and the Piterbarg inequality given in Theorem 8.1 in \cite{Pit2001} 
\BQN \label{eq:Pit}
\mathbb{P}\left(\sup_{(t,\vk{v})\in (t,\vk{v})\in \Delta_x(u)\times (\mathcal{S}_{n-1}/ C_u)}Y_{\vk{b}}(t,\vk{v})>u\right)\leq \eE{\mathbb{Q}} u^{\frac{2(n+1)}{\gamma\wedge 2}-1}\exp\left(-\frac{u^2}{2(1-(\mmu )^2\mathbb{Q}_1)^2}\right) 
\EQN
holds for all $u$ large, with some positive constants $\mathbb{Q}, \mathbb{Q}_1$.
In the following we shall give the asymptotics of $\pi(u)$, from which we shall see that the right-hand side  of   \eqref{eq:Pit} is asymptotically negligible. Thus we conclude that
\BQNY
\mathbb{P}\left(\sup_{(t,\vk{v})\in\mathcal{G}_x(u)}Y_{\vk{b}}(t,\vk{v})>u\right)=\pi(u)(1+o(1)),~~~u\rw\IF.
\EQNY
To this end, we partition $C_u$ into sets of small diameters. To make it more precise, we resort to  the following 
polar coordinates, i.e., for any $\vk{v}\in\mathcal{S}_{n-1}$
\BQNY
v_n= \sin(\theta_1), v_{n-1}=\cos(\theta_1)\sin(\theta_2) \ldot  v_{k+1}=\sin(\theta_{n-k})\prod_{i=1}^{n-k-1} \cos(\theta_i)
\ldot v_1=\prod_{i=1}^{n-1} \cos(\theta_i),
\EQNY
with $\theta_i\in[-{\pi}/{2}, {\pi}/{2}], 1\leq i\leq n-2, \theta_{n-1}\in[0,2\pi).$
We divide $[-{\pi}/{2},{\pi}/{2}]^{k-2}\times [0, 2\pi]$ into several small cubes with length of edge  $R>0,$ denoted by $\{B_j\}_{j\in\NN}$, with
 \BQNY
 \NN=  \{j\in \mathbb{Z}: B_j\subset [-{\pi}/{2},{\pi}/{2}]^{k-2}\times [0, 2\pi]\}.
 \EQNY
Therefore, the corresponding partition of $C_{u}$ can be represented as $\{D_j\}_{j\in \NN}$ with
$$D_j=\Biggl\{\vk{v}=\vk{v}(\theta): (\theta_{n-k+1},\cdots,\theta_{n-1})\in B_j, v_i\in[-\mmu, \mmu], k+1\leq i\leq n
\Biggr\},$$
where
$$\vk{v}(\theta)=(\cos(\theta_1)\cdots \cos(\theta_{n-2} )\cos(\theta_{n-1}),\cdots,\sin(\theta_1)) .$$
Further, set
 \BQNY
 \NN_1=  \{j\in \mathbb{Z}: B_j\subset [- {\pi}/{2}+R,{\pi}/{2}-R]^{k-2}\times [R, 2\pi-R]\}.
 \EQNY
It follows from Bonferroni's inequality that
\BQN\label{EBON}
\sum_{j\in\NN_1}\mathbb{P}\left(\sup_{\Delta_x(u)\times D_j}Y_{\vk{b}}(t,\vk{v})>u\right)-\Pi(u)\leq\pi(u)\leq \sum_{j\in\NN}\mathbb{P}\left(\sup_{\Delta_x(u)\times D_j}Y_{\vk{b}}(t,\vk{v})>u\right),
\EQN
where
$$\Pi(u)=\sum_{i<j\in\NN_1}\mathbb{P}\left(\sup_{\Delta_x(u)\times D_i}Y_{\vk{b}}(t,\vk{v})>u, \sup_{\Delta_x(u)\times D_j}Y_{\vk{b}}(t,\vk{v})>u\right).$$
Since
$$\Cov\left(Y_{\vk{b}}(t,\vk{v}), Y_{\vk{b}}(t',\vk{v}')\right)=\sum_{i=1}^nb_i v_iv_i'\Cov(X(t), X(t'))$$
 the Gaussian random field $Y_{\vk{b}}(t,\vk{v})$ is rotational invariant in law with respect to $\vk{v}$ under the orthogonal matrix
  \BQN\label{MA}
A=\left(
\begin{array}{ccc}
\tilde{A}_k&O\\
O&E_{n-k}
\end{array}
\right),
\EQN
where $\widetilde{A}_k$ is any $k\times k$ orthogonal matrix  and $E_{n-k}$ is the $(n-k)\times (n-k)$ unit matrix.
Hence, for any $j\in\NN$ there exists a orthogonal matrix $A_j$ of the form \eqref{MA}  such that $(1,0,\cdots,0)\in A_jD_j$, and thus
\BQNY
\mathbb{P}\left(\sup_{\Delta_x(u)\times D_j}Y_{\vk{b}}(t,\vk{v})>u\right)&=&\mathbb{P}\left(\sup_{\Delta_x(u)\times D_j}Y_{\vk{b}}(t,A_j\vk{v})>u\right)\\
&=&\mathbb{P}\left(\sup_{\Delta_x(u)\times A_jD_j}Y_{\vk{b}}(t,\vk{v})>u\right).
\EQNY
 Suppose $(1,0,\cdots,0)\in D_0$. Therefore, for the summand in \eqref{EBON} it is  sufficient to consider the case $j=0$. Define  projections \eE{$g_l: \R^l\rw \R^{l-1}, 1 < l \le n$ with $g_l(v_1 \ldot v_l)=(v_2,\cdots,v_l).$}
For $u$ large enough and $R$ sufficient small, we have
\BQNY
\left(g_k( \sqrt{1-(n-k)( \mmu)^2 } D^k_0)\right)\times \left[-\mmu, \mmu\right]^{n-k}
\subset g_n(D_0)\subset \left(g_k D^k_0\right)\times \left[-\mmu, \mmu\right]^{n-k}:=\kal{E}_u(D^k_0),
\EQNY
where $D^k_0=\{(\cos(\theta_1)\cdots \cos(\theta_{k-1}),\cdots, \cos(\theta_1) \sin(\theta_2),
\sin(\theta_1)), (\theta_1,\cdots,\theta_{k-1})\in B_0\}$ is a subset of the $(k-1)$-dimensional unit sphere. Clearly,
for any $R>0$ small enough there exists $\epsilon_R\in (0,1)$ such that
\BQN\label{EMT}
 && (1-\epsilon_R)\Vol\left(D^k_0\right)\leq \Vol\left(g_kD^k_0\right)\leq (1+\epsilon_R)\Vol \left(D^k_0\right),\nonumber \\
&& \Vol \left( g_k(\sqrt{1-(n-k)( \mmu)^2 } D^k_0)\right)
\to \Vol\left(g_kD^k_0\right)
\EQN
hold as $u\rw\IF$.
Consequently, since the projection  $g_n$ on $D_0$ is one-to-one for $u$ large enough and $R$ sufficient small, we conclude that
\BQNY
\mathbb{P}\left(\sup_{(t,\tilde{\vk{v}})\in \Delta_x^{(1)}(u)}Y_{\vk{b}}(t,\tilde{\vk{v}})>u\right)\leq\mathbb{P}\left(\sup_{(t,\vk{v})\in\Delta_x(u)\times D_0}Y_{\vk{b}}(t,\vk{v})>u\right)\leq\mathbb{P}\left(\sup_{(t,\tilde{\vk{v}})\in\Delta_x^{(2)}(u)}Y_{\vk{b}}(t,\tilde{\vk{v}})>u\right),
\EQNY
where $\tilde{\vk{v}}=(v_2 \ldot v_n)$ and
\BQNY
\Delta_x^{(1)}(u) = \Delta_x(u)\times  \kal{E}_u(\sqrt{1-(n-k)( \mmu)^2 } D^k_0), \quad 
\Delta_x^{(2)}(u)=\Delta_x(u)\times \kal{E}_u(D^k_0).
\EQNY
Define  next independent centered stationary Gaussian process  $\{X_1^{\pm\epsilon}(t),t\ge0\}$
and centered homogeneous (stationary) Gaussian random field  $\{X_2^{\pm\epsilon}(\tilde{\vk{v}}),\tilde{\vk{v}}\in\R^{n-1}\}$ with unit variances and correlation functions satisfying
\BQNY
\rho_1^{\pm\epsilon}(t)=1-(1\pm\epsilon)2d\abs{t}^\nu(1+o(1)), t\rw 0,~~~ \rho_2^{\pm\epsilon}(\tilde{\vk{v}})=1-(1\pm\epsilon)\sum_{i=2}^n b_i^2 v_i^2(1+o(1)), \tilde{\vk{v}}\rw \vk{0}.
\EQNY
Then  $\tilde{X}^{\pm\epsilon}(t,\tilde{\vk{v}})=\frac{X_1^{\pm\epsilon}(t)+X_2^{\pm\epsilon}(\tilde{\vk{v}})}{\sqrt{2}}, t\ge0, \tilde{\vk{v}}\in\R^{n-1},$
is a centered homogeneous Gaussian random field with unit variance and correlation function satisfying
\BQNY
\rho^{\pm\epsilon}(t,\tilde{\vk{v}})=1-(1\pm\epsilon)d\abs{t}^\nu(1+o(1))-(1\pm\epsilon)\sum_{i=2}^n\frac{b_i^2 }{2 }v_i^2(1+o(1)),~~~(t,\tilde{\vk{v}})\rw (t_0,\vk{0}).
\EQNY
Further, we have that
\BQN\label{PC}
r_{Y_{\vk{b}}}(t,\vk{v},t',\vk{v}')=1-d\abs{t-t'}^\nu\oo-\sum_{i=2}^n\frac{b_i ^2}{2}(v_i-v'_i)^2\oo  
\EQN
holds as $t, t'\rw t_0, \vk{v}, \vk{v}'\rw (1,0,\cdots,0).$  This can be established as in Lemma 9 in \cite{PI1994b}.
Therefore, by Slepian lemma (cf. \cite{Pit96})  we derive,  for $u$ sufficient large and $R>0$ small enough
\BQNY
\mathbb{P}\left(\sup_{(t,\tilde{\vk{v}})\in \Delta_x^{(2)}(u)}Y_{\vk{b}}(t,\tilde{\vk{v}})>u\right)&\leq&\mathbb{P}\left(\sup_{(t,\tilde{\vk{v}})\in\Delta_x^{(2)}(u)} \frac{{\tilde{X}^{+\epsilon}}(t,\tilde{\vk{v}})}{(1+(1-\epsilon)a|t-t_0|^\mu)(1+(1-\epsilon)\sum_{i=k+1}^n\frac{1-b_i^2}{2}v_i^2)}>u\right)\\
\mathbb{P}\left(\sup_{(t,\tilde{\vk{v}})\in \Delta_x^{(1)}(u)}Y_{\vk{b}}(t,\tilde{\vk{v}})>u\right)&\geq&\mathbb{P}\left(\sup_{(t,\tilde{\vk{v}})\in\Delta_x^{(1)}(u)} \frac{{\tilde{X}^{-\epsilon}}(t,\tilde{\vk{v}})}{(1+(1+\epsilon)a|t-t_0|^\mu)(1+(1+\epsilon)\sum_{i=k+1}^n\frac{1-b_i^2}{2}v_i^2)}>u\right).\\
\EQNY
Consequently, we have from \nelem{GP} that
 \BQNY
 \mathbb{P}\left(\sup_{(t,\tilde{\vk{v}})\in \Delta_x^{(2)}(u)}Y_{\vk{b}}(t,\tilde{\vk{v}})>u\right)&\le& a(\en) (2\pi)^{-k/2} \prod_{i=k+1}^n (1-b_i^2)^{-1/2} \Vol(g_kD_0^k)\\
  &&\times\HMM_{\nu,\mu,d,a}(x_1,x_2) u^{k-2+(2/\nu-2/\mu)_+}e^{-\frac{u^2}{2}}\oo
\EQNY
holds as $u\to\IF$, where  $a(\epsilon)\rw1$ as $\epsilon\rw0$. Since further $\Vol( S_{k-1})= \sum_{j\in\NN} \Vol(D_j)=2\pi^{k/2}/\Gamma(k/2)$ we conclude that
\BQN\label{I1}
 \limsup_{u\rw\IF}
\frac{ \sum_{j\in\NN} \mathbb{P}\left(\sup_{\Delta_x(u)\times D_j}Y_{\vk{b}}(t,\vk{v})>u\right)}{ \prod_{i=k+1}^n(1-b_i^2)^{-1/2} \HMM_{\nu,\mu,d,a}(x_1,x_2)   u^{(2/\nu-2/\mu)_+} \Upsilon_k(u^2) }\leq a'(\epsilon, \en_R)
\EQN
with $a'(\epsilon,\en_R)\rw1$ as $\epsilon\rw0$ and $R\to 0$. Similarly, we have
\BQN\label{I2}
 \liminf_{u\rw\IF}
\frac{ \sum_{j\in\NN_1} \mathbb{P}\left(\sup_{\Delta_x(u)\times D_j}Y_{\vk{b}}(t,\vk{v})>u\right)}{   \prod_{i=k+1}^n(1-b_i^2)^{-1/2} \HMM_{\nu,\mu,d,a}(x_1,x_2)   u^{(2/\nu-2/\mu)_+} \Upsilon_k(u^2) }\ge b'(\epsilon, \en_R)\frac{\sum_{j\in\NN_1} \Vol(D_j)}{\Vol(\mathcal{S}_{n-1})}
\EQN
with $b'(\epsilon,\en_R)\rw1$ as $\epsilon\rw0$ and $R\to 0$.

Next we show that 
\BQN\label{Piu}
\limsup_{R\to 0}\limsup_{u\to\IF}\frac{\Pi(u)}{ u^{(2/\nu-2/\mu)_+} \Upsilon_k(u^2) }=0.
\EQN
We have for any $u>0$
\BQNY
\Pi(u)\leq\sum_{i<j\in\NN_1, D_i\cap D_j=\emptyset}p_{ij}(u) +
\sum_{i<j\in\NN_1, D_i\cap D_j\not =\emptyset}p_{ij}(u)=:\Pi_1(u)+\Pi_2(u),
\EQNY
where
$$p_{ij}(u):=\mathbb{P}\left(\sup_{(t,\vk{v})\in\Delta_x(u)\times D_i}Y_{\vk{b}}(t,\vk{v})>u, \sup_{(t,\vk{v})\in\Delta_x(u)\times D_j}Y_{\vk{b}}(t,\vk{v})>u\right).$$
We first consider the sum taking over $D_i\cap D_j=\emptyset.$ The following standard upper bound
\BQNY
p_{ij}(u) \le \mathbb{P}\left(\sup_{(t,\vk{v},t',\vk{w})\in\Delta_x(u)\times D_i\times\Delta_x(u)\times D_j} \Bigl( \overline{Y_{\vk{b}}}(t,\vk{v})+\overline{Y_{\vk{b}}}(t',\vk{w} )\Bigr)>2u\right)
\EQNY
will be crucial for the proof, where $ \overline{Y_{\vk{b}}}(t,\vk{v})= {Y_{\vk{b}}(t,\vk{v})}/{\sigma_{ Y_{\vk{b}}}(t,\vk{v})}$.  Since $D_i\cap D_j=\emptyset$, in view of \eqref{PCorr} we have that, for any
$(t,\vk{v})\in\Delta_x(u)\times D_i, (t',\vk{w})\in\Delta_x(u)\times D_j$,
\BQNY
\E{(\overline{Y_{\vk{b}}}(t,\vk{v})+\overline{Y_{\vk{b}}}(s,\vk{w}))^2} = 4-2(1-r_{Y_{\vk{b}}}((t,\vk{v},t',\vk{w}))   
\leq 4(1-\delta_1)
\EQNY
holds with  some $\delta_1\in(0,1)$ independent of $j$. Thus,
in the light of Borell-TIS inequality, there exists a common positive constant $\QQ$ such that, for all $u>\QQ$
\BQNY
\mathbb{P}\left(\sup_{(t,\vk{v},s,\vk{w})\in\Delta_x(u)\times D_i\times\Delta_x(u)\times D_j}
\Bigl(\overline{Y_{\vk{b}}}(t,\vk{v})+\overline{Y_{\vk{b}}}(s,\vk{w})\Bigr)
>2u\right)
\leq e^{-\frac{(2u-\QQ)^2}{8(1-\delta_1)}}
\EQNY
holds for all $D_i, D_j, i,j\in\NN$ satisfying $D_i\cap D_j=\emptyset$. 
Consequently,
\BQN\label{I4}
\Pi_1(u)\leq N_R^2 e^{-\frac{(2u-\QQ)^2}{8(1-\delta_1)}}
\EQN
with $N_R$ representing the number of $i$ in $\NN_1.$
Next, for the other sum taking over $D_i\cap D_j\neq\emptyset$ we consider
 first the special summand when $i=0$, so that $(1,0,\cdots,0)\in D_0\cup D_j$. By using the projection $g_n$, we have
\BQNY
p_{0j}(u)&=&\mathbb{P}\left(\sup_{(t,\widetilde{\vk{v}})\in\Delta_x(u)\times \left(g_n(D_0)\right)}Y_{\vk{b}}(t,\widetilde{\vk{v}})>u,\sup_{(s,\widetilde{w})\in \Delta_x(u)\times \left(g_n(D_j)\right)}Y_{\vk{b}}(s,\widetilde{\vk{w}})>u\right)\\
&\leq&\mathbb{P}\left(\sup_{(t,\widetilde{\vk{v}})\in\Delta_x(u)\times \kal{E}_u(D^k_0) 
}Y_{\vk{b}}(t,\widetilde{\vk{v}})>u\right)
+\mathbb{P}\left(\sup_{(t,\widetilde{\vk{v}})\in\Delta_x(u)\times \kal{E}_u(D^k_0)
}Y_{\vk{b}}(t,\widetilde{\vk{v}})>u\right)\\
&&-\mathbb{P}\left(\sup_{(t,\widetilde{\vk{v}})\in\Delta_x(u)\times \kal{E}_u( \sqrt{1-(n-k)(\mmu)^2} D^k_0\cup D^k_j)
}Y_{\vk{b}}(t,\widetilde{\vk{v}})>u\right)
\EQNY
with $\widetilde{\vk{w}}=(w_2,\cdots,w_n)$, and $D_j^k=\{(\cos(\theta_1)\cdots \cos(\theta_{k-1}),\cdots,
\cos(\theta_1)\sin(\theta_2), \sin(\theta_1)), (\theta_1,\cdots,\theta_{k-1})\in B_j\}.$
Further, with the aid of \nelem{GP} we get
\BQNY
 &&\limsup_{u\rw\IF}\frac{p_{0j}(u)}
{(2\pi)^{-k/2} \prod_{i=k+1}^n (1-b_i^2)^{-1/2}   \HMM_{\nu,\mu,d,a}(x_1,x_2) u^{k-2+(2/\nu-2/\mu)_+}e^{-\frac{u^2}{2}}}\\
&&\leq a(\epsilon)\Vol(g_k(D_0^k))+a(\epsilon)\Vol(g_k(D_j^k))-b(\epsilon)\Vol(g_k(D_0^k\cup D_j^k)),
\EQNY
where $a(\en), b(\en)$ are two positive functions such that $a(\en)\to1, b(\en)\to1$ as $\en\to0$. Since for any fixed $i\in\NN_1$, 
 $\#\{j\in\NN_1, D_j^k\cap D_i^k\neq\emptyset\}$ = \# $\{j\in\NN_1, B_j\cap B_i\neq\emptyset\}$ $\leq 3^{k-1}.$
Thus we conclude that
\BQN\label{I5}
&& \limsup_{u\rw\IF}\frac{\Pi_2(u)}{(2\pi)^{-k/2} \prod_{i=k+1}^n (1-b_i)^{-1/2}   \HMM_{\nu,\mu,d,a}(x_1,x_2) u^{k-2+(2/\nu-2/\mu)_+}e^{-\frac{u^2}{2}}}\nonumber\\
&&\leq  3^{k }(c(\epsilon,\en_R)-c'(\epsilon,\en_R))\Vol(S_{k-1})
\EQN
with $c(\epsilon,\en_R)\to1$ and $c'(\epsilon,\en_R)\to1$ as $\en\to0$ and $ R\to 0$. Therefore, we obtain from \eqref{I4} and \eqref{I5} that \eqref{Piu} holds. Consequently, the claim follows from \eqref{I1}, \eqref{I2} and \eqref{Piu} by letting $\en\to0$ and $R\to0$.
This completes the proof. \QED

\COM{
\prooflem{GF} We follow similar arguments as in the proof of Lemma 12.2.1 in \cite{leadbetter1983extremes}. Write first  $t_i\in[0,1], 1\leq i\leq n,$ in dyadic form as
$$t_i=\sum_{j=1}^\IF a_j^{(i)}2^{-j}, 1\leq i\leq n,$$
where $a_j^{(i)}=0 \ \text{or} \ 1$, and write
\BQNY
\xi(\vk{t})=\xi(\vk{0})+\sum_{p=1}^\IF \left(\xi(\vk{t}^{(p)})-\xi(\vk{t}^{(p-1)})\right),  \quad \vk{t}^{(p)}=(\sum_{j=1}^p a_j^{(1)}2^{-j},\cdots,\sum_{j=1}^p a_j^{(n)}2^{-j}).
\EQNY
Further set
 $$\xi(\vk{k},\vk{a},p)=|\xi(k_12^{-p}+a_12^{-p-1},\cdots,k_n2^{-p}+a_n2^{-p-1})-\xi(k_12^{-p},\cdots,k_n2^{-p})|$$
for $p=0,1,2,\cdots, \vk{k}=(k_1,\cdots,k_n), 0\leq k_i\leq 2^p-1$ and $\vk{a}=(a_1,\cdots, a_n)$ with $a_i=0 \ \text{or} \ 1$.  Since $\xi(\vk{t})$ is normal with mean zero, we have 
\BQN\label{APP2}
\mathbb{P}\left(\frac{\xi(\vk{k},\vk{a},p)}{\left(\mathbb{E}\left(\xi(\vk{k},\vk{a},p)\right)^2\right)^{1/2}}\geq x\right)\leq e^{-x^2/2},\ \ \ \forall x>0.
\EQN
Consequently, for $A=\bigcup_{p=0}^\IF\left\{\max_{0\leq k_i\leq 2^p-1, a_i\in\{0, 1\},1\leq i\leq n}\frac{\xi(\vk{k},\vk{a},p)}{\left(\mathbb{E}\left(\xi(\vk{k},\vk{a},p)\right)^2\right)^{1/2}}\geq (2\beta(p+1))^{1/2}\right\}$ we obtain
\BQNY
\mathbb{P}(A)&\leq&\sum_{p=0}^\IF\sum_{i=1}^n\sum_{k_i=0}^{2^p-1}\sum_{a_i\in\{0,1\}}\mathbb{P}\left(\frac{\xi(\vk{k},\vk{a},p)}{\left(\mathbb{E}\left(\xi(\vk{k},\vk{a},p)\right)^2\right)^{1/2}}\geq (2\beta(p+1))^{1/2}\right)\\
&\leq&\sum_{p=0}^\IF 2^{n(p+1)}e^{-\beta(p+1)}=\frac{2^ne^{-\beta}}{1-e^{n\ln 2-\beta}},
\EQNY
implying thus 
$$\mathbb{P}(A)\leq 2^{n+1}e^{-\beta}$$
 for all $\beta \in (0,\IF)$.
Using  (\ref{APP1}) we have
\BQNY
\mathbb{E}\left(\xi(\vk{k},\vk{a},p)\right)^2\leq \QQ\sum_{i=1}^n2^{-(p+1)\alpha_i}\leq n\QQ2^{-(p+1)\alpha_0},
\quad \alpha_0:= \min_{1 \le i \le n} \alpha_i,
\EQNY
from which we know on the complementary event $A^c$ of $A$
\BQNY
|\xi(\vk{k},\vk{a},p)|\leq \left(n\QQ2^{-(p+1)\alpha_0}\right)^{1/2}(2\beta (p+1))^{1/2}.
\EQNY
Consequently, on $A^c$
\BQNY
|\xi(\vk{t})-\xi(0)|&\leq& \sum_{p=1}^\IF \left|\xi(\vk{t}^{(p)})-\xi(\vk{t}^{(p-1)})\right|\\
&\leq&\sum_{p=0}^\IF\left(n\QQ2^{-(p+1)\alpha_0}\right)^{1/2}(2\beta (p+1))^{1/2}\\
&= &\sqrt{\beta \QQ}  \sum_{p=0}^\IF\left(n(p+1)2^{-(p+1)\alpha_0 +1}\right)^{1/2}
=  \left(\frac{\beta \QQ}{4c^*}\right)^{1/2},
\EQNY
where
$ c^* = \Bigl( 2 n^{1/2} \sum_{p=0}^\IF\left((p+1)2^{-(p+1)\alpha_0 +1}\right)^{1/2}\Bigr)^{-2}.$ Therefore,
\BQNY
\mathbb{P}\left(\sup_{\vk{t}\in[0,1]^n}|\xi(\vk{t})-\xi(0)|>\left(\frac{\beta \QQ}{4c^*}\right)^{1/2}\right)\leq\pk{A}\le 2^{n+1}e^{-\beta}.
\EQNY
Choosing $\beta=c^*x^2/\QQ$ in the above and  using the inequality (\ref{APP2}) we conclude that
\BQNY
\mathbb{P}\left(\sup_{\vk{t}\in[0,1]^n}\xi(t)>x\right)&\leq& \mathbb{P}\left(\sup_{\vk{t}\in[0,1]^n}|\xi(t)-\xi(0)|>x/2\right)+\mathbb{P}\left(\xi(0)>x/2\right)\\
&\leq&2^{n+1}e^{-c^*x^2/\QQ}+2^{-1}e^{-x^2/(8\sigma^2)},
\EQNY
establishing the proof. \QED
}


{\bf Proof of \eqref{eq:Corr}:}
We only present the proof for $r_{Y_{\vk{b}}}(t,\vk{v},s, \vk{u})$, since the proof  for $r_{Z_{\vk{b}}}(t,\vk{v},s, \vk{u})$  follows similarly. In the following, all the asymptotics are meant for $s,t\to t_0$, $\sum_{i=k+1}^n v_i^2\to 0$ and $\sum_{i=k+1}^n u_i^2\to 0$.

Denoting $\vk{v}_1=(v_1,\cdots,v_k,0,\cdots,0)$ and $\vk{u}_1=(u_1,\cdots,u_k,0,\cdots,0)$ we have
\BQN\label{eq:co1}
1-r_{Y_{\vk{b}}}(t,\vk{v},s, \vk{u})&=&\frac{\sum_{i=1}^n\E{\left(b_iX_i(t)v_i-b_iX_i(s)u_i\right)^2}-\left(\sigma_{Y_{\vk{b}}} (t,\vk{v})-\sigma_{Y_{\vk{b}}} (s,\vk{u})\right)^2}{2\sigma_{Y_{\vk{b}}} (t,\vk{v})\sigma_{Y_{\vk{b}}} (s,\vk{u})}\nonumber\\
&=&\frac{\sigma_{Y_{\vk{b}}} (t,\vk{v}_1)\sigma_{Y_{\vk{b}}} (s,\vk{u}_1)}{\sigma_{Y_{\vk{b}}} (t,\vk{v})\sigma_{Y_{\vk{b}}} (s,\vk{u})}\xi_1(t,\vk{v},s, \vk{u}) +\xi_2(t,\vk{v},s, \vk{u})  +\xi_3(t,\vk{v},s, \vk{u}),
\EQN
where
 \BQNY
&& \xi_1(t,\vk{v},s, \vk{u})=\frac{\sum_{i=1}^k\E{\left(X_i(t)v_i-X_i(s)u_i\right)^2}-\left(\sigma_{Y_{\vk{b}}} (t,\vk{v}_1)-\sigma_{Y_{\vk{b}}} (s,\vk{u}_1)\right)^2}{2\sigma_{Y_{\vk{b}}} (t,\vk{v}_1)\sigma_{Y_{\vk{b}}} (s,\vk{u}_1)},\\
&&\xi_2(t,\vk{v},s, \vk{u})=\frac{\sum_{i=k+1}^n\E{\left(X_i(t)v_i-X_i(s)u_i\right)^2}}{2\sigma_{Y_{\vk{b}}} (t,\vk{v})\sigma_{Y_{\vk{b}}} (s,\vk{u})},\\
&&\xi_3(t,\vk{v},s, \vk{u})=\frac{\left(\sigma_{Y_{\vk{b}}} (t,\vk{v}_1)-\sigma_{Y_{\vk{b}}} (s,\vk{u}_1)\right)^2 -\left(\sigma_{Y_{\vk{b}}} (t,\vk{v})-\sigma_{Y_{\vk{b}}} (s,\vk{u})\right)^2}{2\sigma_{Y_{\vk{b}}} (t,\vk{v})\sigma_{Y_{\vk{b}}} (s,\vk{u})}.
\EQNY
By assumption {\bf II} we have
\BQN\label{eq:co2}
 \xi_1(t,\vk{v},s, \vk{u}) &=&1-\frac{\E{X_1(t)X_1(s)}}{\sigma_{X_1}(t)\sigma_{X_1}(s)}\frac{\sum_{i=1}^kv_iu_i}{\sqrt{\sum_{i=1}^kv_i^2}\sqrt{\sum_{i=1}^k}u_i^2}\nonumber\\
&=& d\abs{t-s}^\nu\oo +\frac{1}{2}\sum_{i=1}^k(v_i-u_i)^2\oo+o\left(\sum_{i=k+1}^n(u_i-v_i)^2\right).
\EQN
Further, we have
\BQNY\label{eq:AA}
\sum_{i=k+1}^n\E{\left(b_iX_i(t)v_i-b_iX_i(s)u_i\right)^2} 
&=&\sum_{i=k+1}^nb_i^2\Big(\E{\left(X_i(t)-X_i(s)\right)^2}v_i^2+\E{\left(X_i(s)v_i-X_i(s)u_i\right)^2} \nonumber\\
&& +2\E{(X_i(t)v_i-X_i(s)v_i)(X_i(s)v_i-X_i(s)u_i)}\Big).
\EQNY
Next we deal with the last three terms on the right-hand side in turn. By assumption {\bf III} and the fact that $\gamma_1\geq \nu$ we get
\BQNY
\sum_{i=k+1}^nb_i^2\E{\left(X_i(t)-X_i(s)\right)^2}v_i^2\leq G_1\sum_{i=k+1}^nb_i^2v_i^2|t-s|^{\gamma_1}=o(|t-s|^{\nu}),
\EQNY
 and
\BQNY
\sum_{i=k+1}^nb_i^2\E{\left(X_i(s)v_i-X_i(s)u_i\right)^2}=\sum_{i=k+1}^nb_i^2\sigma_{X_{k+1}}^2(s)(u_i-v_i)^2=\sum_{i=k+1}^nb_i^2(u_i-v_i)^2\oo.
\EQNY
In addition, we have
\BQNY
&&2\sum_{i=k+1}^nb_i^2\E{(X_i(t)v_i-X_i(s)v_i)(X_i(s)v_i-X_i(s)u_i)}\\
&&\leq 2\sum_{i=k+1}^nb_i^2\left(\E{(X_i(t)v_i-X_i(s)v_i)^2}\right)^{1/2}\left(\E{(X_i(s)v_i-X_i(s)u_i)^2}\right)^{1/2}\\
&&\leq2\sum_{i=k+1}^nb_i^2 G_1^{1/2}|t-s|^{\gamma_1/2}|v_i||v_i-u_i|\\
&&\le \sum_{i=k+1}^nb_i^2 G_1^{1/2}(|t-s|^{\gamma_1 }|v_i|+\abs{v_i}(v_i-u_i)^2)\\
&&=o(|t-s|^\nu)+o\LT(\sum_{i=k+1}^n(v_i-u_i)^2\RT).
\EQNY
Therefore, we obtain
\BQN\label{eq:co3}
\xi_2(t,\vk{v},s, \vk{u}) 
=\frac{1}{2}\sum_{i=k+1}^nb_i^2(u_i-v_i)^2\oo+o(|t-s|^\nu).
\EQN
Moreover, it follows that
\BQNY
\left(\sigma_{Y_{\vk{b}}} (t,\vk{v}_1)-\sigma_{Y_{\vk{b}}} (s,\vk{u}_1)\right)^2&=&
\left(\sigma_{Y_{\vk{b}}} (t,\vk{v}_1)-\sigma_{Y_{\vk{b}}} (t,\vk{u}_1)\right)^2+\left(\sigma_{Y_{\vk{b}}} (t,\vk{u}_1)-\sigma_{Y_{\vk{b}}}(s,\vk{u}_1)\right)^2 \\
&&+2(\sigma_{Y_{\vk{b}}} (t,\vk{v}_1)-\sigma_{Y_{\vk{b}}} (t,\vk{u}_1))(\sigma_{Y_{\vk{b}}} (t,\vk{u}_1)-\sigma_{Y_{\vk{b}}}(s,\vk{u}_1)).
\EQNY
Direct calculation yields that
\BQNY
\left(\sigma_{Y_{\vk{b}}} (t,\vk{v}_1)-\sigma_{Y_{\vk{b}}} (t,\vk{u}_1)\right)^2&\le&\left(\left(\sum_{i=1}^kv_i^2\right)^{1/2}-\left(\sum_{i=1}^ku_i^2\right)^{1/2}\right)^2\\
&=&a(\vk{v},\vk{u})\sum_{i=k+1}^n(v_i-u_i)^2
\EQNY
with $a(\vk{v},\vk{u})\to 0,$ and
\BQNY
\left(\sigma_{Y_{\vk{b}}} (t,\vk{u}_1)-\sigma_{Y_{\vk{b}}}(s,\vk{u}_1)\right)^2= \left(\sigma_{X_1}(t)-\sigma_{X_1}(s)\right)^2(1+o(1)) 
\EQNY
hold. 
Since further  by assumption {\bf III}
\BQNY
\left(\sigma_{X_1}(t)-\sigma_{X_1}(s)\right)^2\leq G|t-s|^\gamma
\EQNY
we conclude,  by the fact that $\gamma\geq \nu$,
\BQNY
&&2(\sigma_{Y_{\vk{b}}} (t,\vk{v}_1)-\sigma_{Y_{\vk{b}}} (t,\vk{u}_1))(\sigma_{Y_{\vk{b}}} (t,\vk{u}_1)-\sigma_{Y_{\vk{b}}}(s,\vk{u}_1))\\
&&\le (a(\vk{v},\vk{u}))^{-1/2}\left(\sigma_{Y_{\vk{b}}} (t,\vk{v}_1)-\sigma_{Y_{\vk{b}}} (t,\vk{u}_1)\right)^2+(a(\vk{v},\vk{u}))^{1/2}\left(\sigma_{Y_{\vk{b}}} (t,\vk{u}_1)-\sigma_{Y_{\vk{b}}}(s,\vk{u}_1)\right)^2 \\
&&=o\left(|t-s|^\nu\RT)+o\LT(\sum_{i=k+1}^n(v_i-u_i)^2\right)
\EQNY
implying that
\BQNY
\left(\sigma_{Y_{\vk{b}}} (t,\vk{v}_1)-\sigma_{Y_{\vk{b}}} (s,\vk{u}_1)\right)^2=\left(\sigma_{X_1}(t)-\sigma_{X_1}(s)\right)^2(1+o(1))+o\left(|t-s|^\nu\RT)+o\LT(\sum_{i=k+1}^n(v_i-u_i)^2\right).
\EQNY
Similarly,
\BQNY
\left(\sigma_{Y_{\vk{b}}} (t,\vk{v})-\sigma_{Y_{\vk{b}}} (s,\vk{u})\right)^2=\left(\sigma_{X_1}(t)-\sigma_{X_1}(s)\right)^2(1+o(1))+o\left(|t-s|^\nu\RT)+o\LT(\sum_{i=k+1}^n(v_i-u_i)^2\right).
\EQNY
Therefore,
\BQN\label{eq:co4}
\xi_3(t,\vk{v},s, \vk{u})&=&o(\left(\sigma_{X_1}(t)-\sigma_{X_1}(s)\right)^2)+o\left(|t-s|^\nu\RT)+o\LT(\sum_{i=k+1}^n(v_i-u_i)^2\right)\nonumber\\
&=&o\left(|t-s|^\nu\RT)+o\LT(\sum_{i=k+1}^n(v_i-u_i)^2\right).
\EQN
Consequently, the claim in \eqref{eq:Corr} follows by combining (\ref{eq:co1})--(\ref{eq:co4}). 
\QED

\bigskip

{\bf Acknowledgement}: Both authors kindly acknowledge partial
support from Swiss National Science Foundation Project 200021-140633/1,
and
the project RARE -318984   (an FP7  Marie Curie IRSES Fellowship). P. Liu also acknowledges partial support by the National Science Foundation of China 11271385.
\bibliographystyle{plain}

 \bibliography{gausbib}

\end{document}